\documentclass[11pt]{article}
\usepackage[a4paper,margin=27mm]{geometry}
\usepackage[T1]{fontenc}
\usepackage{lmodern}
\usepackage{amsmath,amssymb,amsthm,mathtools}
\usepackage{mathrsfs}
\usepackage{microtype}
\usepackage[colorlinks=true,linkcolor=blue,citecolor=blue,urlcolor=blue]{hyperref}
\usepackage{xcolor}
\newtheorem{lemma}{Lemma}
\newtheorem{theorem}[lemma]{Theorem}
\newtheorem{corollary}[lemma]{Corollary}
\newtheorem{remark}[lemma]{Remark}
\newtheorem{definition}[lemma]{Definition}
\newcommand{\R}{\mathbb R}
\newcommand{\llbracket}{\mathopen{[\![}}
\newcommand{\rrbracket}{\mathclose{]\!]}}
\newcommand{\mres}{\mathbin{\llcorner}}

\allowdisplaybreaks[1]

\title{$C^{1,\alpha}$ Regularity of $W^{2,2}$ $Q$-Valued Maps with Small Laplacian}
\author{Mattia Luchese}
\date{}
\begin{document}
\maketitle
\vspace{-1.5em}

\begin{abstract}
We establish a regularity criterion for $Q$-valued maps $u\in W^{2,2}_{\mathrm{loc}}(\Omega,\mathcal A_Q(\mathbb R^k))$, where $\Omega\subset\mathbb R^n$ is open. If $|\Delta u|\leq\delta_*|D^2u|$ almost everywhere, with $\delta_*=\delta_*(n,k,Q)>0$ sufficiently small, then $u\in C^{1,\alpha}_{\mathrm{loc}}$ for some $\alpha=\alpha(n,k,Q)>0$, with quantitative estimates. The result holds for arbitrary $n\geq2$, $k\geq1$, and $Q\geq1$.

The proof establishes uniform energy decay using a rigidity principle for positive-semidefinite matrix-valued measures with an aperture bound and approximately atomic divergence, extending the target-stress argument for the homogeneity gap in \cite{Luc26}. The argument relies only on first variation identities.

As consequences, we obtain $C^{1,\alpha}$ regularity for Dirichlet-stationary $W^{2,2}$ maps and for locally Lipschitz $W^{2,2}$ maps whose graph varifolds are stationary and have sufficiently small slope. Combining the latter result with Bellettini's Lipschitz sheeting theorem \cite{Bel25}, we establish $C^{1,\alpha}$ sheeting in the multiple-valued sense for two-sided stable minimal immersions whose singular set has vanishing $2$-capacity relative to the immersion, under suitable flatness and mass bounds. Our results allow us to reduce the proof of compactness within this class, under locally uniform mass bounds, to the conjecture that branch sets of stationary $C^{1,\alpha}$ multiple-valued graphs in codimension one have vanishing $2$-capacity.
\end{abstract}

\tableofcontents

\section{Introduction}

Almgren's theory of Dirichlet-minimizing $Q$-valued functions \cite{Alm00}, developed as part of the regularity theory for area-minimizing currents and revisited by De Lellis and Spadaro \cite{DLS11}, establishes local H\"older continuity and a bound of $n-2$ for the Hausdorff dimension of the singular set.
Branching is compatible with energy minimization, and these conclusions do not in general extend to continuity of the first derivative: the $2$-valued square root already provides an example.
For maps which are just stationary with respect to the Dirichlet energy, much less is known. On $2$-dimensional domains, Lin \cite{Lin14} proved continuity, and Hirsch and Spolaor \cite{HS24} subsequently established H\"older continuity and a sharp lower bound of $1$ for the codimension of the singular set.
These results concern maps in the natural energy space $W^{1,2}$.

The present paper studies the additional regularity forced by a second-order Sobolev assumption.
We prove that every Dirichlet-stationary $Q$-valued map in $W^{2,2}_{\mathrm{loc}}$ is locally $C^{1,\alpha}$, in arbitrary domain dimension, target dimension, and multiplicity.
More generally, for every $n\geq2$ and $k,Q\geq1$, Theorem~\ref{thm:small-trace-regularity} gives constants $\delta_*=\delta_*(n,k,Q)>0$, $\alpha=\alpha(n,k,Q)\in(0,1)$ and $C=C(n,k,Q)<\infty$ such that every map
\[
u\in W^{2,2}\bigl(B_1^n,\mathcal A_Q(\R^k)\bigr)
\]
satisfying
\[
|\Delta u|\leq\delta_*|D^2u|
\qquad\text{almost everywhere in }B_1^n
\]
admits a $C^{1,\alpha}_{\mathrm{loc}}$ representative, with
\[
\|u\|_{C^{1,\alpha}(B_{1/2}^n)}
\leq
C\|u\|_{W^{2,2}(B_1^n)}.
\]
We stress that no further assumptions on $u$ are required.
In particular, $u$ is not assumed to satisfy any prescribed differential equation or variational identity.

To the best of the author's knowledge, this criterion does not appear to have been stated previously in this precise form even for single-valued maps.
In that case, however, it admits a much simpler proof by harmonic replacement and the standard decay estimates for the Hessian energy of harmonic maps; see \cite[Lemma~1.I and the proof of Lemma~2.I]{Giu67} for the corresponding argument in the scalar case.

This result complements the Schauder theory for multiple-valued elliptic systems.
Simon and Wickramasekera \cite{SW16} established Schauder estimates for $2$-valued symmetric solutions of linear systems and proved the optimal $C^{1,1/2}$ regularity of stationary $2$-valued graphs initially assumed to be $C^{1,\alpha}$.
Building on their method, the author established in \cite{Luc26} a priori $C^{1,\alpha}$ and $C^{2,\alpha}$ estimates for $Q$-valued quasilinear systems which are perturbations of the Laplace system, with arbitrary multiplicity and in arbitrary dimension and codimension.
Those estimates concern solutions already belonging to the corresponding H\"older classes.
The present paper obtains the initial $C^{1,\alpha}$ regularity from Sobolev assumptions, under the small-Laplacian condition above.
In the Dirichlet-stationary case, it also gives a quantitative converse to the implication established through the Bochner inequality in \cite{Luc26}, namely that $C^1$ maps belong locally to $W^{2,2}$.

The connection with \cite{Luc26} also concerns the method of proof.
We use a generalisation of the variational identities established there for the differential of $W^{2,2}$ $Q$-valued maps.
The basic idea behind the energy-decay argument is already present in the target-stress argument used there to establish the homogeneity gap, although in that setting it is formulated for maps already known to be homogeneous.
That argument relies on a rigidity lemma for positive-semidefinite matrix-valued measures with an aperture bound and atomic divergence, showing that the divergence must vanish.
Here we extend this rigidity lemma to approximately atomic divergence, allowing the atomic loads to move and come together, and use it to establish uniform energy decay.
The argument relies only on first variation identities and does not use minimisation properties, either directly or through previous regularity results.

The small Laplacian formulation also applies to the minimal surface system.
For a locally Lipschitz $W^{2,2}$ map whose graph varifold is stationary, the approximate second jets satisfy the minimal surface system in non-divergence form almost everywhere, giving
\[
|\Delta u|\leq \|Du\|_{L^\infty}^2|D^2u|.
\]
Consequently, sufficiently small slope implies the hypothesis of the regularity criterion.
This yields $C^{1,\alpha}$ regularity for stationary multiple-valued graphs under the $W^{2,2}$ assumption, with estimates depending only on the dimensions and multiplicity.
The resulting graphs belong to the class treated by the Schauder estimates and the branch-set dimension theorem of \cite{Luc26}.

A geometric application arises in the sheeting theory for stable minimal hypersurfaces.
The theorem of Schoen and Simon \cite{SS81} gives smooth sheeting in the embedded setting.
For immersed hypersurfaces, branching must also be allowed.
Wickramasekera \cite{Wic08} established $2$-valued $C^{1,\alpha}$ sheeting near a multiplicity-$2$ hyperplane, under an assumption of locally finite $\mathcal H^{n-2}$-measure on the singular set.
Recently, Bellettini \cite{Bel25} obtained a Lipschitz multiple-valued graphical description near a hyperplane with arbitrary multiplicity, for two-sided stable minimal immersions whose singular set has vanishing $2$-capacity relative to the immersion.
In particular, singular sets of locally finite $\mathcal H^{n-2}$-measure are permitted.
His argument uses intrinsic differential inequalities and De Giorgi iteration.

Even more recently, Wang and Zhang \cite{WZ26} have developed this geometric analysis further by establishing branched sheeting near stationary classical cones and unions of hyperplanes whose common intersection has dimension at least $n-2$.
Their theorem gives stationary Lipschitz multiple-valued graphs over the constituent hyperplanes and leads to regularity and precompactness results controlling the non-branch singular set.
In the discussion following Remark~5.1 of \cite{Bel25}, Bellettini explicitly raises the question of whether the Lipschitz conclusion can be improved to $C^{1,\alpha}$.

Corollary~\ref{cor:stable-sheeting} gives a positive answer to this question.
The connection with the analytic criterion is provided by a capacitary extension argument: after Lipschitz sheeting, the relative capacity assumption implies that the projected singular set has vanishing Sobolev $2$-capacity in the base.
This allows the first jet of the graphical representation to be extended in $W^{1,2}$, giving the required $W^{2,2}$ regularity of the map.
The small-slope result then yields $C^{1,\alpha}$ sheeting with quantitative estimates.
The conclusion is understood in the multiple-valued sense and permits both branching and crossings.

In Section~\ref{sec:stable-immersions}, we study the regularity of varifold limits of two-sided stable minimal immersions whose singular sets (in the immersed sense) have vanishing relative $2$-capacity.
Under locally uniform mass bounds, we show that each limit is again represented by a smooth two-sided stable minimal immersion away from an immersed singular set of Hausdorff dimension at most $n-2$.
This dimension bound follows by combining the sheeting theory of \cite{WZ26}, the present regularity estimates and the branch-set dimension bound established in \cite[Theorem~5.7]{Luc26}, whose proof in turn uses the stratification theory of Krummel, Minter and Wickramasekera \cite{KMW26}.
We also prove compactness within the original class under the stronger conjecture that branch sets of stationary $C^{1,\alpha}$ multiple-valued graphs in codimension one have vanishing Sobolev $2$-capacity.

Finally, we point out that, at the same time, independently and using different methods, Sidney Stanbury obtained overlapping regularity results in \cite{Sta26}, relying on Dirichlet minimisation and comparison rather than the first-variation methods used here. In particular, Corollaries~\ref{cor:W22-regularity}, \ref{thm:area-regularity} and \ref{cor:stable-sheeting}, and the singular-set dimension bound in Theorem~\ref{thm:conditional-compactness}, correspond to Proposition~2.4, Theorems~1.14--1.15 and Corollary~1.16 in his paper.

\subsection*{Acknowledgements}

This work was supported by the UK Engineering and Physical Sciences Research Council (EPSRC) Studentship EP/T517847/1 and by the Cambridge Trust. It was carried out while the author was a PhD student at the University of Cambridge.

The author is very grateful to his supervisor, Neshan Wickramasekera, for suggesting a version of the problem that motivated this work, and for his guidance and several helpful discussions on this paper and on related topics. The author also thanks Costante Bellettini for reading an early draft of the paper and making useful suggestions, including the application to compactness for stable minimal immersions, which is developed in the final section.

\section{Regularity of maps with small Laplacian}
\label{sec:dirichlet-core}

All norms and contractions are Euclidean.
For matrices $A,B$ of the same size, we write $A:B:=\operatorname{tr}(A^*B)$, where $A^*$ denotes the Euclidean adjoint of $A$.
Whenever convenient, we identify Euclidean spaces with their duals through the inner product, and in particular we use the same notation for the differential of a function and its gradient.
We denote by $\operatorname{Sym}^+(\R^m)$ the cone of symmetric positive-semidefinite matrices.
For a symmetric matrix-valued measure $\mathsf S$ on $\R^m$, the divergence is understood distributionally:
\[
\langle\operatorname{div}\mathsf S,\Phi\rangle:=-\int_{\R^m}D\Phi:d\mathsf S
\qquad\text{for every }\Phi\in C_c^1(\R^m;\R^m).
\]

\subsection{Rigidity for matrix-valued measures}

\begin{lemma}[Approximately atomic divergence]
\label{lem:moving-atomic-loads}
Let $m,q\geq1$ be integers and let $c\in(0,1)$.
Let $\mathsf S_j\in
\mathcal M\bigl(\R^m;\operatorname{Sym}^+(\R^m)\bigr)\subset \mathcal D'\bigl(\R^m;\R^{m\times m}\bigr)$
be a sequence of finite measures satisfying
\begin{equation}
\label{eq:moving-aperture}
\begin{cases}
    \sup_j\operatorname{tr}\mathsf S_j(\R^m)<\infty,\\
    e^T\mathsf S_je\leq c\,\operatorname{tr}\mathsf S_j
    \qquad\text{as measures }\forall e\in \mathbb S^{m-1},\\
    \operatorname{div}\mathsf S_j
    =\sum_{\ell=1}^q b_{j,\ell}\,\delta_{p_{j,\ell}}+R_j,
\end{cases}
\end{equation}
where the points $p_{j,\ell}\in\R^m$ satisfy
$\sup_{j,\ell}|p_{j,\ell}|<\infty$, $b_{j,\ell}\in\R^m$ are arbitrary, 
and $R_j\in (C^1_c(\R^m, \R^m))^*\subset \mathcal D'(\R^m, \R^m)$ are small in the sense that for every compact set $K\subset\R^m$
\begin{equation}
\label{eq:moving-test-defect}
\|R_j\|_{C^1(K)^*}
:=
\sup_{\substack{
\Phi\in C_c^1(\R^m;\R^m)\\
\operatorname{supp}\Phi\subset K,\ \|\Phi\|_{C^1}\leq1}}
|\langle R_j,\Phi\rangle|
\longrightarrow0,
\end{equation}
where $\|\Phi\|_{C^1}:=\|\Phi\|_\infty+\|D\Phi\|_\infty$.
Then
\begin{equation}
\label{eq:moving-local-zero}
\mathsf S_j\stackrel{*}{\rightharpoonup}0
\qquad\text{as matrix-valued Radon measures on }\R^m.
\end{equation}
\end{lemma}

\begin{proof}
We argue by induction on $q$.

\medskip
\noindent
\emph{Step 1: atomicity of weak limits.}
By the mass bound, after passing to a subsequence we may assume $p_{j,\ell}\longrightarrow p_\ell$, $\mathsf S_j\stackrel{*}{\rightharpoonup}\mathsf S$ locally in $\R^m$.
The limit $\mathsf S$ is a finite positive-semidefinite matrix-valued measure and inherits the aperture bound.
Let $A:=\{a^1,\ldots,a^N\}$ be the set of distinct points among the $p_\ell$.

If $\Phi\in C_c^1(\R^m\setminus A;\R^m)$, then $\Phi(p_{j,\ell})=0$ for every $\ell$ and all sufficiently large $j$.
Hence \eqref{eq:moving-aperture} gives $\int_{\R^m}D\Phi:d\mathsf S_j =-\langle R_j,\Phi\rangle$.
Passing to the limit and using \eqref{eq:moving-test-defect}, we obtain $\int_{\R^m}D\Phi:d\mathsf S=0$, so that
\[
\operatorname{div}\mathsf S=0
\qquad\text{in }\R^m\setminus A.
\]
The distribution $\operatorname{div}\mathsf S$ has order at most one and is supported in the finite set $A$.
Indeed, if a test field $\Phi$ satisfies $\Phi(a^\beta)=D\Phi(a^\beta)=0$ for every $\beta$, multiplying it by cutoffs supported in shrinking neighbourhoods of $A$ leaves the action of $\operatorname{div}\mathsf S$ unchanged and produces fields converging to zero in $C^1$.
Hence $\langle\operatorname{div}\mathsf S,\Phi\rangle=0$.
Its action therefore depends linearly only on the first jets $(\Phi(a^\beta),D\Phi(a^\beta))$.
Since these form a finite-dimensional vector space, there exist $b_\beta\in\R^m$ and $C_\beta\in\R^{m\times m}$ such that
\begin{equation}
\label{eq:moving-point-distribution}
\operatorname{div}\mathsf S
=\sum_{\beta=1}^N
\bigl(b_\beta\delta_{a^\beta}
+\operatorname{div}(C_\beta\delta_{a^\beta})\bigr).
\end{equation}

We identify $C_\beta$ with $\mathsf S(\{a^\beta\})$.
Fix $\Phi\in C_c^1(\R^m;\R^m)$ with $\Phi(0)=0$ and set $\Phi_\rho(y):= \rho\,\Phi\!\left(\frac{y-a^\beta}{\rho}\right)$.
For all sufficiently small $\rho$, its support contains no other point of $A$, so \eqref{eq:moving-point-distribution} gives
\[
-C_\beta:D\Phi(0)
=\langle\operatorname{div}\mathsf S,\Phi_\rho\rangle
=-\int_{\R^m}D\Phi_\rho:d\mathsf S
\xrightarrow[\rho\downarrow0]{}
-\mathsf S(\{a^\beta\}):D\Phi(0),
\]
where the limit follows from dominated convergence, as $D\Phi_\rho$ is uniformly bounded, converges pointwise to zero away from $a^\beta$, and equals $D\Phi(0)$ at $a^\beta$.
Since $D\Phi(0)$ is arbitrary, comparison gives
\[
C_\beta=\mathsf S(\{a^\beta\}).
\]

Thus the derivative-of-Dirac terms are precisely the divergence of the
atomic part of $\mathsf S$. We remove that part by defining
\[
\mathsf S^0:=\mathsf S-\sum_{\beta=1}^N \mathsf S(\{a^\beta\})\,\delta_{a^\beta} =\mathsf S\llcorner(\R^m\setminus A).
\]
The measure $\mathsf S^0$ is positive semidefinite, has no atom at any point
of $A$, and inherits the aperture bound in \eqref{eq:moving-aperture}. Moreover, $\operatorname{div}\mathsf S^0 =\sum_{\beta=1}^Nb_\beta\delta_{a^\beta}$.
Applying \cite[Lemma~3.13]{Luc26} in a ball containing only $a^\beta$,
after translation and rescaling, with load $-b_\beta$, yields $b_\beta=0$. Repeating this at every point of $A$,
we obtain
$$
\operatorname{div}\mathsf S^0=0 \qquad \text{in $\R^m$}.
$$

Since $\mathsf S^0$ is finite, choose $Y_R(y)=\zeta(y/R)y$, where $\zeta=1$ on $B_1^m$ and $\zeta=0$ outside $B_2^m$.
Then $0=\int_{\R^m}DY_R:d\mathsf S^0$.
The matrices $DY_R$ are uniformly bounded and converge pointwise to the
identity. Dominated convergence gives $0=\operatorname{tr}\mathsf S^0(\R^m)$, and hence
$$
\mathsf S^0=0.
$$
Therefore every fixed-scale limit has the form
\begin{equation}
\label{eq:moving-fixed-scale-atomic}
\mathsf S=\sum_{\beta=1}^NC_\beta\delta_{a^\beta}.
\end{equation}

\medskip
\noindent
\emph{Step 2: base case $q=1$.}
Let $p_{j,1}\to a$ and let
$\mathsf S_j\stackrel{*}{\rightharpoonup}\mathsf S$ locally. By
\eqref{eq:moving-fixed-scale-atomic}, $\mathsf S=C_1\delta_a$ for some $C_1\in\operatorname{Sym}^+(\R^m)$. Choose
$\zeta\in C_c^1(\R^m)$ equal to one near $a$. Set $\Phi_j(y):=\zeta(y)(y-p_{j,1})$.
Then $\Phi_j(p_{j,1})=0$.
The fields $\Phi_j$ have uniformly bounded $C^1$ norms and are supported in the fixed compact set $\operatorname{supp}\zeta$.
By \eqref{eq:moving-aperture} and
\eqref{eq:moving-test-defect},
\[
\operatorname{tr}C_1
=\lim_{j\to\infty}\int_{\R^m}D\Phi_j:d\mathsf S_j
=-\lim_{j\to\infty}\langle R_j,\Phi_j\rangle
=0.
\]
Since $C_1$ is positive semidefinite, $C_1=0$. Thus every local limit is zero, proving
\eqref{eq:moving-local-zero} when $q=1$.

\medskip
\noindent
\emph{Step 3: inductive step, $q\ge2$ and $N\ge2$.}
Assume the conclusion proved for all multiplicities smaller than $q$ and
consider a subsequence as in Step 1. Suppose that $N\geq2$. Define $I_\beta:=\{\ell:p_\ell=a^\beta\}$, $q_\beta:=|I_\beta|<q$.
Choose $\chi_\beta\in C_c^1(\R^m)$ such that $0\leq\chi_\beta\leq1$,
$\chi_\beta=1$ near $a^\beta$, and $\chi_\beta=0$ near every
$a^\gamma$ with $\gamma\neq\beta$. Put
\[
\mathsf S_j^\beta:=\chi_\beta\mathsf S_j.
\]
These measures preserve positivity, the mass bound, and the aperture
condition. For all sufficiently large $j$, we have
$\chi_\beta(p_{j,\ell})=1$ for $\ell\in I_\beta$ and
$\chi_\beta(p_{j,\ell})=0$ for $\ell\notin I_\beta$.
The product rule gives
\[
\operatorname{div}\mathsf S_j^\beta
=\chi_\beta\operatorname{div}\mathsf S_j
+\mathsf S_jD\chi_\beta
=\sum_{\ell\in I_\beta}b_{j,\ell}\,\delta_{p_{j,\ell}}
+R_j^\beta,
\]
where $R_j^\beta:=\chi_\beta R_j+\mathsf S_jD\chi_\beta$.
The support of $D\chi_\beta$ is a compact set disjoint from $A$, while the fixed-scale limit is supported in $A$ by \eqref{eq:moving-fixed-scale-atomic}.
Consequently, $\operatorname{tr}\mathsf S_j(\operatorname{supp}D\chi_\beta)\longrightarrow0$.
For every compact set $K\subset\R^m$, multiplication by $\chi_\beta$ preserves support in $K$, and therefore
\[
\|R_j^\beta\|_{C^1(K)^*}
\leq
\|\chi_\beta\|_{C^1}\|R_j\|_{C^1(K)^*}
+\|D\chi_\beta\|_\infty\operatorname{tr}\mathsf S_j(\operatorname{supp}D\chi_\beta)
\longrightarrow0.
\]
Hence the induction hypothesis applies
to $\mathsf S_j^\beta$ with multiplicity $q_\beta$, and $\mathsf S_j^\beta\stackrel{*}{\rightharpoonup}0$.
On the other hand, $\mathsf S_j^\beta\stackrel{*}{\rightharpoonup}
C_\beta\delta_{a^\beta}$.
Thus $C_\beta=0$. Repeating this for every $\beta$ gives
$\mathsf S=0$. Therefore every local limit is zero whenever the limiting
configuration contains at least two distinct points.

\medskip
\noindent
\emph{Step 4: inductive step, $q\ge2$ and $N=1$.}
Assume now that $p_{j,\ell}\longrightarrow a$ for every $\ell$.
By Step 1, every local limit has the form $C_1\delta_a$. We show that $C_1=0$.
If all the $p_{j,\ell}$ coincide along a subsequence, the argument of Step 2
applies directly. We may therefore assume $\rho_j:=\max_\ell|p_{j,\ell}-p_{j,1}|>0$.
Set $c_j:=p_{j,1}$, so that $\rho_j\to0$.

Suppose by contradiction that $C_1\neq0$. Choose a cutoff
$\chi\in C_c^1(\R^m)$ satisfying $0\leq\chi\leq1$ and
\[
\chi=1\text{ on }B_1^m(a),
\qquad
\chi=0\text{ outside }B_{2}^m(a),
\]
and a function $\vartheta\in C^1([0,\infty))$ such that
\[
\vartheta=0\text{ on }[0,1],
\qquad
\vartheta=1\text{ on }[2,\infty),
\qquad
\vartheta'\geq0.
\]
Define
\[
\Phi_j(y):=
\chi(y)\,
\vartheta\!\left(\frac{|y-c_j|}{\rho_j}\right)
(y-c_j).
\]
All points $p_{j,\ell}$ satisfy $|p_{j,\ell}-c_j|\leq\rho_j$, so
$\Phi_j(p_{j,\ell})=0$ for every $\ell$. The $C^1$ norms are uniformly bounded:
on the transition annulus the factor $|y-c_j|/\rho_j$ is bounded.
Moreover, all these fields are supported in the fixed compact set $\operatorname{supp}\chi$, so \eqref{eq:moving-test-defect} applies.
Consequently,
\begin{equation}
\label{eq:moving-collapse-integral-zero}
\int_{\R^m}D\Phi_j:d\mathsf S_j
=-\langle R_j,\Phi_j\rangle
\longrightarrow0.
\end{equation}

Since $\vartheta'\geq0$, the radial field satisfies
$
D\!\left[
\vartheta\!\left(\frac{|y-c_j|}{\rho_j}\right)(y-c_j)
\right]
\geq
\vartheta\!\left(\frac{|y-c_j|}{\rho_j}\right)I.
$
The term in $D\Phi_j$ containing $D\chi$ is uniformly bounded and supported
in a compact set disjoint from $a$, whose $\operatorname{tr}\mathsf S_j$-mass tends to zero.
Positivity and \eqref{eq:moving-collapse-integral-zero} therefore give
\[
0\leq
\operatorname{tr}\mathsf S_j\bigl(B_1^m(a)\setminus B_{2\rho_j}^m(c_j)\bigr)
\leq
\left|
\int_{\R^m}D\Phi_j:d\mathsf S_j
\right|
+C\operatorname{tr}\mathsf S_j(\operatorname{supp}D\chi)
\longrightarrow0.
\]
Since $B_{2\rho_j}^m(c_j)\subset B_1^m(a)$ for all sufficiently large $j$ and
$\operatorname{tr}\mathsf S_j(B_1^m(a))\to\operatorname{tr}C_1$, we conclude that
\[
\operatorname{tr}\mathsf S_j\bigl(B_{2\rho_j}^m(c_j)\bigr)
\longrightarrow\operatorname{tr}C_1>0.
\]

Rescale the target around $c_j$ at scale $\rho_j$. Define
\[
\eta_{c_j,\rho_j}(y):=\frac{y-c_j}{\rho_j},
\qquad
\widetilde{\mathsf S}_j:=(\eta_{c_j,\rho_j})_\#\mathsf S_j,
\qquad
\widetilde p_{j,\ell}:=\eta_{c_j,\rho_j}(p_{j,\ell}).
\]
Then $\widetilde p_{j,1}=0$ and $\max_\ell|\widetilde p_{j,\ell}|=1$. After a subsequence, at least two
of the limiting points $\widetilde p_\ell$ are distinct. The mass bound and the aperture
condition are unchanged by this push-forward.

The rescaled measures satisfy
\[
\operatorname{div}\widetilde{\mathsf S}_j
=\sum_{\ell=1}^q\rho_j b_{j,\ell}\,
  \delta_{\widetilde p_{j,\ell}}+\widetilde R_j,
\qquad
\widetilde R_j:=\rho_j(\eta_{c_j,\rho_j})_\#R_j.
\]
Fix a compact set $K\subset\R^m$.
Since $c_j\to a$ and $\rho_j\to0$, we have
\[
c_j+\rho_jK\subset\overline{B_1^m(a)}
\]
for all sufficiently large $j$.
Thus, for every $\Phi\in C_c^1(\R^m;\R^m)$ supported in $K$ with $\|\Phi\|_{C^1}\leq1$,
\[
\begin{aligned}
|\langle\widetilde R_j,\Phi\rangle|
&=\bigl|\langle R_j,\rho_j\Phi\circ\eta_{c_j,\rho_j}\rangle\bigr|\\
&\leq\|R_j\|_{C^1(\overline{B_1^m(a)})^*}
\bigl(\rho_j\|\Phi\|_\infty+\|D\Phi\|_\infty\bigr)\\
&\leq\|R_j\|_{C^1(\overline{B_1^m(a)})^*},
\end{aligned}
\]
since $D[\rho_j\Phi\circ\eta_{c_j,\rho_j}]=D\Phi\circ\eta_{c_j,\rho_j}$ and $\rho_j\leq1$ eventually.
Taking the supremum over these tests gives
\[
\|\widetilde R_j\|_{C^1(K)^*}
\leq\|R_j\|_{C^1(\overline{B_1^m(a)})^*}
\longrightarrow0.
\]
Thus all the hypotheses of Step~3 are satisfied, and its conclusion
contradicts
\[
\operatorname{tr}\widetilde{\mathsf S}_j(B_2^m)
=
\operatorname{tr}\mathsf S_j\bigl(B_{2\rho_j}^m(c_j)\bigr)
\longrightarrow\operatorname{tr}C_1>0.
\]
A compactly supported continuous nonnegative test equal to one on
$\overline{B_2^m}$ makes this contradiction immediate.
Hence $C_1=0$. The induction is complete, and
\eqref{eq:moving-local-zero} follows.
\end{proof}

\subsection{The regularity theorem}

We use the Sobolev spaces and first-jet conventions of \cite[Section~2]{Luc26}.
Let $n\geq2$, $k,Q\geq1$, and let $\Omega\subset\R^n$ be open.
Set $E:=(\R^k)^n$, naturally identified with $\R^{k\times n}$.
A map $u\in W^{1,2}_{\mathrm{loc}}(\Omega,\mathcal A_Q(\R^k))$ belongs to $W^{2,2}_{\mathrm{loc}}$ if its approximate first jet satisfies $j_1u:=\sum_{\ell=1}^Q\llbracket(u^\ell,Du^\ell)\rrbracket\in W^{1,2}_{\mathrm{loc}}(\Omega,\mathcal A_Q(\R^k\times E))$.
For such a map, write $G:=Du\in W^{1,2}_{\mathrm{loc}}(\Omega,\mathcal A_Q(E))$, with $G_i^\ell=D_i u^\ell\in\R^k$.
We regard $DG^\ell=(D_jG_i^\ell)_{i,j=1}^n=(D_{ij}u^\ell)_{i,j=1}^n$ as an $n\times n$ matrix with entries in $\R^k$, representing the linear map $DG^\ell\in\operatorname{Hom}(\R^n,E)$ given by $(DG^\ell v)_i=\sum_{j=1}^n D_{ij}u^\ell\,v_j$.
Here $(DG^\ell)^*$ denotes its Euclidean adjoint.
Thus $DG^\ell(DG^\ell)^*\in\operatorname{Hom}(E,E) \simeq \R^{k\times k}\otimes \R^{n\times n}$ is an $n\times n$ block matrix with entries in $\R^{k\times k}$,
given by $\bigl[DG^\ell(DG^\ell)^*\bigr]_{ij}
=\sum_{p=1}^n D_{ip}u^\ell \, ( D_{jp}u^\ell)^T \in \R^{k\times k}$, where we regard vectors in $\R^k$ as column vectors.
Similarly, for $\chi\in C^1(\Omega)$, $(DG^\ell D\chi)_i =\sum_{j=1}^n D_{ij}u^\ell\,D_j\chi \in \R^k$, so that $DG^\ell D\chi\in E$.
We write $\Delta u^\ell:=\sum_{i=1}^n D_{ii}u^\ell$, $\Delta u:=\sum_{\ell=1}^Q\llbracket\Delta u^\ell\rrbracket$, where the trace is taken in the domain indices, and set $|\Delta u|^2:=\sum_\ell|\Delta u^\ell|^2$, $|D^2u|^2:=\sum_\ell|D^2u^\ell|^2=|DG|^2$.
For $Y\in E$, let $Y_i\in\R^k$ denote its $i$th column, and write $\operatorname{div}\Psi:=\sum_iD_i\Psi_i$ for an $E$-valued field.
For a symmetric tensor $T=(T_{ij})$, we set $(\operatorname{div}T)_j:=\sum_iD_iT_{ij}$.
In the identities below, inner products between vectors in $\R^k$ are denoted by juxtaposition, and derivatives of $\Psi(x,G^\ell(x))$ are taken with respect to $x$ after composition.
Dirichlet stationarity means both outer and inner stationarity.

\medskip
The following lemma extends the computations of \cite[Theorem~3.11]{Luc26}, where the condition $\Delta u=0$ was used to prove that $G=Du$ is Dirichlet-stationary.
The aim is to quantify the failure of outer and inner stationarity of $G$ in terms of $\Delta u$, with errors that will be controlled below under the assumption $|\Delta u|\leq\delta|D^2u|$.
No new ideas are needed: we repeat the same computations as in \cite[Theorem~3.11]{Luc26}, keeping track of the terms that vanished when $\Delta u=0$.
The algebraic computations yielding the explicit error terms were initially carried out using ChatGPT.

\begin{lemma}[Variation identities for the differential]
\label{lem:weak-identities}
Let $u\in W^{2,2}_{\mathrm{loc}}(\Omega,\mathcal A_Q(\R^k))$ and set $G:=Du$.
The identity
\begin{equation}
\label{eq:weak-hodge}
\int_\Omega\sum_\ell DG^\ell:D[\Psi(x,G^\ell)]
=\int_\Omega\sum_\ell\Delta u^\ell\operatorname{div}[\Psi(x,G^\ell)]
\end{equation}
holds for every $\Psi\in C_c^1(\Omega\times E;E)$.
The same identity holds for $\Psi\in C^1(\Omega\times E;E)$ supported in $K\times E$ for some compact set $K\Subset\Omega$, provided that $|D_Y\Psi(x,Y)|\leq C$, $|\Psi(x,Y)|+|D_x\Psi(x,Y)|\leq C(1+|Y|)$
for some constant $C<\infty$.

Setting $T_{ij}(G):=|DG|^2\delta_{ij}-2\sum_\ell D_iG^\ell\cdot D_jG^\ell$, $B_{ij}:=|\Delta u|^2\delta_{ij}-2\sum_\ell\Delta u^\ell D_iG_j^\ell$,
we have
\begin{equation}
\label{eq:rellich}
\int_\Omega T(G):DX=\int_\Omega B:DX
\end{equation}
for every $X\in C_c^1(\Omega;\R^n)$.

Moreover, if $B_R^n\Subset\Omega$, then $r\mapsto D_G(r):=r^{2-n}\int_{B_r^n}|DG|^2\,dx$ is locally absolutely continuous on $(0,R)$ and, for almost every $r\in(0,R)$,
\begin{equation}
\label{eq:radial-identity}
\begin{aligned}
D_G'(r)
={}&2r^{2-n}\int_{\partial B_r^n}|\partial_rG|^2\\
&-(n-2)r^{1-n}\int_{B_r^n}|\Delta u|^2
+r^{2-n}\int_{\partial B_r^n}\sum_\ell
\bigl[|\Delta u^\ell|^2-2\Delta u^\ell D^2u^\ell[\nu,\nu]\bigr],
\end{aligned}
\end{equation}
where $\nu:=x/|x|$ and $\partial_rG^\ell:=DG^\ell[\nu]$.
\end{lemma}

\begin{proof}
The componentwise symmetry of the second differential follows from \cite[Lemma~3.10]{Luc26}.
Applying the paired Stokes formula \cite[Lemma~3.9]{Luc26} componentwise gives
\[
\begin{aligned}
\int_\Omega\sum_{\ell,i,j}D_jG_i^\ell D_j[\Psi_i(x,G^\ell)]
&=\int_\Omega\sum_{\ell,i,j}D_iG_j^\ell D_j[\Psi_i(x,G^\ell)]\\
&=\int_\Omega\sum_{\ell,i,j}D_jG_j^\ell D_i[\Psi_i(x,G^\ell)]\\
&=\int_\Omega\sum_\ell\Delta u^\ell\operatorname{div}[\Psi(x,G^\ell)],
\end{aligned}
\]
which proves \eqref{eq:weak-hodge} for compactly supported fields.

For fields satisfying the stated growth bounds, choose $\vartheta\in C_c^\infty(E)$ with $\vartheta=1$ on $\{|Y|\leq1\}$ and $\vartheta=0$ on $\{|Y|\geq2\}$, and apply \eqref{eq:weak-hodge} to $\Psi_M(x,Y):=\vartheta(Y/M)\Psi(x,Y)$.
The growth bounds imply, uniformly for $M\geq1$, $\bigl|D[\Psi_M(x,G^\ell)]\bigr|
\leq C\bigl(1+|G^\ell|+|DG^\ell|\bigr)$.
Since $G,DG,\Delta u\in L^2(K)$, dominated convergence allows passage to the limit as $M\to\infty$, proving the extension of \eqref{eq:weak-hodge}.

For the inner variation identity, define $T_{ij}(u):=|G|^2\delta_{ij}-2\sum_\ell G_i^\ell G_j^\ell$.
For readability, we also write $\Delta u\,G_j$ for $\sum_\ell\Delta u^\ell G_j^\ell$.
The product rule and the symmetry $D_iG_j^\ell=D_jG_i^\ell$ give, in distributions,
\begin{equation}
\label{eq:inner-var-u-divergence}
\sum_iD_iT_{ij}(u)=-2\Delta u\,G_j.
\end{equation}
Apply \eqref{eq:weak-hodge} with $\Psi(x,Y)=\varphi(x)Y$ and, for fixed $i,j$, with the field whose $r$th column is $\Psi_r(x,Y)=\tfrac12\varphi(x)\bigl(\delta_{rj}Y_i+\delta_{ri}Y_j\bigr)$, $\varphi\in C_c^\infty(\Omega)$.
Expanding these identities as in \cite[Theorem~3.11]{Luc26}, but without assuming $\Delta u=0$, gives
\begin{equation}
\label{eq:inner-var-G}
T_{ij}(G)
=\tfrac12\Delta T_{ij}(u)+B_{ij}
-\delta_{ij}\sum_rD_r(\Delta u\,G_r)
+D_i(\Delta u\,G_j)+D_j(\Delta u\,G_i)
\end{equation}
in distributions.
All undifferentiated products are locally integrable.
Applying $\sum_iD_i$, the last three terms reduce to $\Delta(\Delta u\,G_j)$, which equals $-\tfrac12\sum_iD_i\Delta T_{ij}(u)$ by \eqref{eq:inner-var-u-divergence}.
Thus $\operatorname{div}[T(G)]=\operatorname{div}B$ in $\mathcal D'(\Omega;\R^n)$, proving \eqref{eq:rellich}.

The local absolute continuity of $D_G$ follows from coarea.
Choose $\eta\in C^\infty(\R)$ with $0\leq\eta\leq1$, $\eta=1$ on $(-\infty,0]$, and $\eta=0$ on $[1,\infty)$.
Testing \eqref{eq:rellich} with $X_\varepsilon(x):=\eta((|x|-r)/\varepsilon)x$, where $0<\varepsilon<R-r$, and letting $\varepsilon\downarrow0$ as in the proof of \cite[Proposition~3.2]{DLS11}, we obtain, for almost every $r\in(0,R)$,
\[
\int_{B_r^n}\operatorname{tr}\bigl(T(G)-B\bigr)
=r\int_{\partial B_r^n}\bigl(T(G)-B\bigr)[\nu,\nu].
\]
Substituting the definitions of $T(G)$ and $B$ and differentiating $D_G(r)$ yields \eqref{eq:radial-identity}.
\end{proof}

\begin{theorem}[$W^{2,2}$-to-$C^{1,\alpha}$ regularity]
\label{thm:small-trace-regularity}
Let $n\geq2$ and $k,Q\geq1$.
There exist $\delta_*>0$, $\alpha\in(0,1)$ and $C<\infty$, depending only on $n,k,Q$, with the following property.
Every map $u\in W^{2,2}(B_1^n,\mathcal A_Q(\R^k))$ satisfying $|\Delta u|\leq\delta_*|D^2u|$ almost everywhere in $B_1^n$ admits a representative in $C^{1,\alpha}_{\mathrm{loc}}(B_1^n)$, with
\begin{equation}\label{eq:small-trace-holder}
[Du]_{C^{0,\alpha}(B_{1/2}^n)}\leq C\|D^2u\|_{L^2(B_1^n)}.
\end{equation}
As a consequence,
\begin{equation}\label{eq:small-trace-full-norm}
\|u\|_{C^{1,\alpha}(B_{1/2}^n)}\leq C\|u\|_{W^{2,2}(B_1^n)}.
\end{equation}
\end{theorem}

\begin{proof}
Throughout this proof, the Euclidean scalar product in $\mathbb R^k$ is left implicit.
Assume $|\Delta u|\leq\delta|D^2u|$, with $0<\delta\leq1/2$ to be chosen below.
Set $G:=Du$, $E:=(\R^k)^n$, and $D_G(r):=r^{2-n}\int_{B_r^n}|DG|^2\,dx$ for $0<r\leq1$.
We may assume $D_G(1)>0$, since otherwise $u$ is affine and the conclusion is immediate.

In Steps~1--3 we establish a uniform contraction estimate for the rescaled Hessian energy.
In Step~4 we iterate this estimate to obtain the regularity conclusion and the quantitative bounds.

\medskip
\noindent\emph{Step 1: from outer variation to the target stress.}
The weak Hodge identity \eqref{eq:weak-hodge} gives, for every $\Psi\in C_c^1(B_1^n\times E;E)$,
\[
\int_{B_1^n}\sum_{\ell=1}^Q DG^\ell:D[\Psi(\cdot,G^\ell)]\,dx
=\int_{B_1^n}\sum_{\ell=1}^Q\Delta u^\ell\operatorname{div}[\Psi(\cdot,G^\ell)]\,dx,
\]
where the derivatives are taken with respect to $x$ after composition.
Fix a radial $\chi\in C_c^\infty(B_1^n)$ with $0\leq\chi\leq1$ and $\chi=1$ on $B_{1/2}^n$.
Taking $\Psi(x,Y)=\chi(x)\Phi(Y)$, with $\Phi\in C_c^1(E;E)$, and applying the chain rule, we obtain
\[
\begin{aligned}
\int_{B_1^n}\chi\sum_{\ell=1}^Q D\Phi(G^\ell):DG^\ell(DG^\ell)^*\,dx
&=-\int_{B_1^n}\sum_{\ell=1}^Q\Phi(G^\ell)\cdot DG^\ell D\chi\,dx\\
&\quad+\int_{B_1^n}\sum_{\ell=1}^Q\Delta u^\ell\operatorname{div}[\chi\Phi(G^\ell)]\,dx.
\end{aligned}
\]
For each sheet and target component, define
\[
H_0^{\ell,a}:=D^2u^{\ell,a}-\frac{\Delta u^{\ell,a}}n I_n,
\]
and write $H_0^\ell\in\operatorname{Hom}(\mathbb R^n,E)$ for the tensor with these components.
Replacing $DG^\ell(DG^\ell)^*$ by $H_0^\ell(H_0^\ell)^*$ gives
\[
\begin{aligned}
\int_{B_1^n}\chi\sum_{\ell=1}^Q D\Phi(G^\ell):H_0^\ell(H_0^\ell)^*\,dx
&=-\int_{B_1^n}\sum_{\ell=1}^Q\Phi(G^\ell)\cdot DG^\ell D\chi\,dx\\
&\quad+\int_{B_1^n}\sum_{\ell=1}^Q\Delta u^\ell\operatorname{div}[\chi\Phi(G^\ell)]\,dx\\
&\quad+\int_{B_1^n}\chi\sum_{\ell=1}^Q D\Phi(G^\ell):\bigl(H_0^\ell(H_0^\ell)^*-DG^\ell(DG^\ell)^*\bigr)\,dx.
\end{aligned}
\]
Define the normalized target stress
\[
\mathsf S:=\frac1{D_G(1)}\int_{B_1^n}\chi\sum_{\ell=1}^Q H_0^\ell(H_0^\ell)^*\,\delta_{G^\ell}\,dx.
\]
The preceding identity then becomes
\[
\begin{aligned}
\langle\operatorname{div}\mathsf S,\Phi\rangle
&=-\int_E D\Phi:d\mathsf S\\
&=\frac1{D_G(1)}\int_{B_1^n}\sum_{\ell=1}^Q\Phi(G^\ell)\cdot DG^\ell D\chi\,dx\\
&\quad-\frac1{D_G(1)}\int_{B_1^n}\sum_{\ell=1}^Q\Delta u^\ell\operatorname{div}[\chi\Phi(G^\ell)]\,dx\\
&\quad+\frac1{D_G(1)}\int_{B_1^n}\chi\sum_{\ell=1}^Q D\Phi(G^\ell):\bigl(DG^\ell(DG^\ell)^*-H_0^\ell(H_0^\ell)^*\bigr)\,dx.
\end{aligned}
\]

Since the component matrices $H_0^{\ell,a}$ are symmetric and trace-free, \cite[Lemma~3.15]{Luc26} gives
\[
H_0^\ell(H_0^\ell)^*\leq\frac{n-1}{n}|H_0^\ell|^2 I_E.
\]
Thus $\mathsf S$ is positive semidefinite and has aperture at most $(n-1)/n<1$.
Moreover, the orthogonal trace decomposition gives
\[
\sum_{\ell=1}^Q|H_0^\ell|^2=|DG|^2-\frac1n|\Delta u|^2,
\]
and therefore
\[
\operatorname{tr}\mathsf S(E)
=\frac1{D_G(1)}\int_{B_1^n}\chi\sum_{\ell=1}^Q|H_0^\ell|^2\,dx
=\frac1{D_G(1)}\int_{B_1^n}\chi\left(|DG|^2-\frac1n|\Delta u|^2\right)\,dx.
\]
Using $0\leq\chi\leq1$, $\chi=1$ on $B_{1/2}^n$, and $|\Delta u|\leq\delta|DG|$, we obtain
\begin{equation}
\label{eq:energy-stress-mass}
\left(1-\frac{\delta^2}{n}\right)2^{2-n}\frac{D_G(1/2)}{D_G(1)}\leq\operatorname{tr}\mathsf S(E)\leq1.
\end{equation}

\medskip
\noindent\emph{Step 2: estimate of the remainder.}
We first decompose the divergence of $\mathsf S$ into an atomic part plus a remainder.
Choose
\[
P=\sum_{\ell=1}^Q\llbracket p_\ell\rrbracket
\in\operatorname*{argmin}_{T\in\mathcal A_Q(E)}\int_{B_1^n}\mathcal G(G,T)^2.
\]
Choose measurable labels realizing an optimal matching of $G$ to $P$, and relabel the corresponding approximate differentials accordingly, so that
\[
\sum_{\ell=1}^Q|G^\ell-p_\ell|^2=\mathcal G(G,P)^2
\qquad\text{a.e. on }B_1^n.
\]
Set
\[
b_\ell:=\frac1{D_G(1)}\int_{B_1^n}\bigl(DG^\ell D\chi-\Delta u^\ell\otimes D\chi\bigr)\,dx,
\]
and define $R$ by
\[
\operatorname{div}\mathsf S=\sum_{\ell=1}^Qb_\ell\,\delta_{p_\ell}+R.
\]
Thus, for every $\Phi\in C_c^1(E;E)$,
\begin{equation}
\label{eq:stress-remainder-identity}
\begin{aligned}
\langle R,\Phi\rangle
={}&\frac1{D_G(1)}\int_{B_1^n}\sum_{\ell=1}^Q
\bigl(\Phi(G^\ell)-\Phi(p_\ell)\bigr)\cdot
\bigl(DG^\ell D\chi-\Delta u^\ell\otimes D\chi\bigr)\,dx\\
&-\frac1{D_G(1)}\int_{B_1^n}\chi\sum_{\ell=1}^Q
\Delta u^\ell\operatorname{div}[\Phi(G^\ell)]\,dx\\
&+\frac1{D_G(1)}\int_{B_1^n}\chi\sum_{\ell=1}^Q
D\Phi(G^\ell):\bigl(DG^\ell(DG^\ell)^*-H_0^\ell(H_0^\ell)^*\bigr)\,dx.
\end{aligned}
\end{equation}
The orthogonal trace decomposition gives
\begin{equation}
\label{eq:stress-trace-error}
\sum_{\ell=1}^Q\bigl|DG^\ell(DG^\ell)^*-H_0^\ell(H_0^\ell)^*\bigr|
\leq C|\Delta u||DG|
\leq C\delta|DG|^2.
\end{equation}
Moreover, the chain rule and Cauchy--Schwarz give
\[
\sum_{\ell=1}^Q\bigl|\Delta u^\ell\operatorname{div}[\Phi(G^\ell)]\bigr|
\leq C\|D\Phi\|_\infty|\Delta u||DG|
\leq C\delta\|D\Phi\|_\infty|DG|^2.
\]
Using this bound and \eqref{eq:stress-trace-error} in \eqref{eq:stress-remainder-identity}, we obtain
\begin{equation}
\label{eq:stress-defect}
\begin{aligned}
|\langle R,\Phi\rangle|
&\leq\frac1{D_G(1)}
\left|\int_{B_1^n}\sum_{\ell=1}^Q
\bigl(\Phi(G^\ell)-\Phi(p_\ell)\bigr)\cdot
\bigl(DG^\ell D\chi-\Delta u^\ell\otimes D\chi\bigr)\,dx\right|\\
&\quad+C\delta\|D\Phi\|_\infty.
\end{aligned}
\end{equation}

The choice of $P$ and the metric Poincar\'e inequality give
\begin{equation}
\label{eq:energy-poincare}
\int_{B_1^n}\mathcal G(G,P)^2
\leq C\int_{B_1^n}|DG|^2
=C D_G(1).
\end{equation}
Moreover, optimal matching gives
\[
\sum_{\ell=1}^Q|\Phi(G^\ell)-\Phi(p_\ell)|^2
\leq\|D\Phi\|_\infty^2\mathcal G(G,P)^2.
\]
Since $\chi$ is radial and $D\chi$ is supported in $B_1^n\setminus B_{1/2}^n$, we also have
\[
\int_{B_1^n}\sum_{\ell=1}^Q\bigl|DG^\ell D\chi-\Delta u^\ell\otimes D\chi\bigr|^2
\leq C\left[\int_{B_1^n\setminus B_{1/2}^n}|\partial_rG|^2+\delta^2D_G(1)\right].
\]
Thus Cauchy--Schwarz in \eqref{eq:stress-defect} yields
\begin{equation}
\label{eq:stress-remainder-radial}
|\langle R,\Phi\rangle|
\leq C\|D\Phi\|_\infty
\left[
\left(\frac{\displaystyle\int_{B_1^n\setminus B_{1/2}^n}|\partial_rG|^2}{D_G(1)}\right)^{1/2}
+\delta
\right].
\end{equation}

We now relate the radial derivative to the energy drop.
Integrating \eqref{eq:radial-identity} over $[1/2,s]$, letting $s\uparrow1$, and using the small-Laplacian bound yields
\begin{equation}
\label{eq:radial-integrated}
\left|D_G(1)-D_G(1/2)-2\int_{B_1^n\setminus B_{1/2}^n}|x|^{2-n}|\partial_rG|^2\right|
\leq C\delta D_G(1).
\end{equation}
In particular, when $\Delta u=0$, the weighted radial energy is exactly half the energy drop.
In general, since $n\geq2$,
\begin{equation}
\label{eq:energy-small-radial}
\frac{\displaystyle\int_{B_1^n\setminus B_{1/2}^n}|\partial_rG|^2}{D_G(1)}
\leq C\left[\left(1-\frac{D_G(1/2)}{D_G(1)}\right)_++\delta\right],
\end{equation}
where $t_+:=\max\{t,0\}$.
Combining this with \eqref{eq:stress-remainder-radial}, and using $\delta\leq1/2$, we obtain
\begin{equation}
\label{eq:stress-remainder-gap}
\|R\|_{(C_0^1)^*}
\leq C\left[\left(1-\frac{D_G(1/2)}{D_G(1)}\right)_++\delta\right]^{1/2},
\end{equation}
where $\|R\|_{(C_0^1)^*}:=\sup\{|\langle R,\Phi\rangle|:\Phi\in C_c^1(E;E),\ \|\Phi\|_\infty+\|D\Phi\|_\infty\leq1\}$.
Thus a small relative energy drop and a small Laplacian bound force the remainder to be small.

We shall also use \eqref{eq:stress-defect} and \eqref{eq:stress-remainder-radial} for fields of the form $\Phi(Y)=Y-\Phi_c(Y)$, with $\Phi_c\in C_c^1(E;E)$.
These fields have bounded derivatives and at most linear growth, so the extension in Lemma~\ref{lem:weak-identities} applies to $\Psi(x,Y)=\chi(x)\Phi(Y)$.
The preceding identities and estimates therefore remain valid for these fields.

\medskip
\noindent\emph{Step 3: uniform energy decay.}
We now prove that there exist $\delta_*=\delta_*(n,k,Q)\in(0,1/2]$ and $\theta=\theta(n,k,Q)\in(0,1)$, such that
\begin{equation}
\label{eq:small-trace-contraction}
D_G(1/2)\leq\theta D_G(1)
\end{equation}
whenever $|\Delta u|\leq\delta_*|DG|$ almost everywhere in $B_1^n$.
Suppose otherwise.
Then there are maps $u_j$ on $B_1^n$, with $G_j=Du_j$ and $\delta_j\downarrow0$, such that
\[
|\Delta u_j|\leq\delta_j|DG_j|,
\qquad
D_{G_j}(1)>0,
\qquad
\frac{D_{G_j}(1/2)}{D_{G_j}(1)}>1-j^{-1}.
\]
Multiplying each $u_j$ by a positive constant, we may normalize
\begin{equation}
\label{eq:energy-sequence-normalization}
\int_{B_1^n}|G_j|^2=1.
\end{equation}
This normalization will keep the reference configurations bounded, while \eqref{eq:energy-stress-mass} will give a positive lower bound for the stress mass.

Apply the constructions of Steps~1 and~2 to each $u_j$, and denote the resulting objects by $P_j=\sum_{\ell=1}^Q \llbracket p_{j,\ell} \rrbracket$, $\mathsf S_j$, $b_{j,\ell}$ and $R_j$ such that
\[
\operatorname{div} \mathsf S_j = \sum_{\ell=1}^Q b_{j,\ell} \delta_{p_{j,\ell}} + R_j.
\]
Comparing $P_j$ with $Q\llbracket0\rrbracket$ gives
\[
\int_{B_1^n}\mathcal G(G_j,P_j)^2\leq1.
\]
The triangle inequality then yields
\[
|B_1^n||P_j|^2
\leq2\int_{B_1^n}|G_j|^2+2\int_{B_1^n}\mathcal G(G_j,P_j)^2
\leq4.
\]
Thus the points $p_{j,\ell}$ are uniformly bounded.

By \eqref{eq:stress-remainder-gap},
\[
\|R_j\|_{(C_0^1)^*}
\leq C(j^{-1}+\delta_j)^{1/2}\longrightarrow0.
\]
The stresses are positive semidefinite, have uniformly bounded mass and have aperture at most $(n-1)/n<1$.
Hence Lemma~\ref{lem:moving-atomic-loads} (identifying $E$ with $\mathbb R^m$, where $m=kn$) gives
\[
\mathsf S_j\rightharpoonup^*0
\qquad\text{locally on }E.
\]
On the other hand, \eqref{eq:energy-stress-mass} gives
\begin{equation}
\label{eq:stress-positive-mass}
\liminf_{j\to\infty}\operatorname{tr}\mathsf S_j(E)\geq2^{2-n}>0.
\end{equation}

To get a contradiction, it remains to exclude escape of this mass to infinity.
As the $p_{j,\ell}$ are uniformly bounded, choose $\zeta\in C_c^\infty(E)$ equal to one on a ball containing every $p_{j,\ell}$, and set
\[
\Phi(Y):=Y-\zeta(Y)Y.
\]
This field vanishes at every $p_{j,\ell}$ and has bounded derivative.
Since the estimates of Step~2 remain valid for this field, \eqref{eq:stress-remainder-radial} and \eqref{eq:energy-small-radial} give
\[
\left|\int_E D\Phi:d\mathsf S_j\right|
\leq C\|D\Phi\|_\infty(j^{-1}+\delta_j)^{1/2}
\longrightarrow0.
\]
Local convergence also gives
\[
\int_E D(\zeta(Y)Y):d\mathsf S_j\longrightarrow0.
\]
Since $D\Phi+D(\zeta(Y)Y)=I_E$, we conclude that
\[
\operatorname{tr}\mathsf S_j(E)
=\int_E D\Phi:d\mathsf S_j+\int_E D(\zeta(Y)Y):d\mathsf S_j
\longrightarrow0,
\]
contradicting \eqref{eq:stress-positive-mass}.
This proves \eqref{eq:small-trace-contraction}.

\medskip\noindent
\emph{Step 4: regularity and quantitative estimates.}
We now take $\delta\leq\delta_*$, with $\delta_*$ as in Step~3.
Enlarge $\theta$ if necessary so that $1/2\leq\theta<1$, and set
\[
\alpha:=-\frac{\log\theta}{2\log2}\in(0,1/2].
\]
For every ball $B_R^n(x)\Subset B_1^n$, the rescaled map $v(z):=R^{-1}u(x+Rz)$ satisfies the same small-Laplacian condition.
Applying \eqref{eq:small-trace-contraction} to these rescalings and iterating gives
\[
\int_{B_{2^{-m}R}^n(x)}|DG|^2\leq2^{-m(n-2+2\alpha)}\int_{B_R^n(x)}|DG|^2,
\qquad m=0,1,\ldots.
\]
For $2^{-(m+1)}R<s\leq2^{-m}R$, comparison on nested balls yields
\begin{equation}\label{eq:small-trace-decay}
s^{2-n}\int_{B_s^n(x)}|DG|^2\leq C\left(\frac{s}{R}\right)^{2\alpha}R^{2-n}\int_{B_R^n(x)}|DG|^2,
\qquad 0<s\leq R,\quad B_R^n(x)\Subset B_1^n.
\end{equation}
Here one may take $C=2^{n-2+2\alpha}$.

Compose $G$ with an Almgren bi-Lipschitz embedding and apply the Campanato--Morrey theorem to \eqref{eq:small-trace-decay}.
This gives a locally $\alpha$-H\"older representative of the embedded map.
Since the image of the embedding is closed, this representative takes values in that image everywhere and can be pulled back.
We thus obtain a locally $\alpha$-H\"older representative of $G$ on $B_1^n$.
In particular, taking $R=1/2$ in \eqref{eq:small-trace-decay} for every centre $x\in B_{1/2}^n$ gives
\[
[G]_{C^{0,\alpha}(B_{1/2}^n)}\leq C\|DG\|_{L^2(B_1^n)},
\]
which is estimate \eqref{eq:small-trace-holder}.

We next identify this continuous differential with the Taylor first jet of $u$.
By \cite[Definition~2.8 and Remark~2.9(ii)]{Luc26}, $j_1u\in W^{1,2}_{\mathrm{loc}}$ and $|D(j_1u)|^2=|G|^2+|DG|^2$.
Since $G$ is locally bounded and $\alpha<1$, this identity and \eqref{eq:small-trace-decay} give the same local Morrey bound for $j_1u$.
An Almgren embedding therefore gives a locally $\alpha$-H\"older representative of the paired jet, and hence of $u$.
Its gradient projection agrees everywhere with the continuous representative of $G$.

Sobolev slicing and the one-dimensional selection theorem \cite[Proposition~1.2]{DLS11}, applied to $j_1u$, give absolutely continuous selections $(u^\ell,G^\ell)$ on almost every line in each fixed direction $e$.
Compatibility of the approximate jet and Sobolev locality give
\[
(u^\ell)'=G^\ell[e]\qquad\text{a.e.}
\]
Fix $x_0\in B_1^n$ and a neighbourhood $U\Subset B_1^n$ of $x_0$.
For $h\neq0$ sufficiently small, set $e=h/|h|$ and choose good segments $[x_q,x_q+h]\subset U$ with $x_q\to x_0$.
Along each segment, the image of every selected gradient is contained in the union of the $Q$ closed balls of radius $[G]_{C^{0,\alpha}(U)}|h|^\alpha$ centred at the initial gradient values.
By continuity, this image lies in a single connected component of the union, whose diameter is at most $2Q[G]_{C^{0,\alpha}(U)}|h|^\alpha$.
Integrating the selection identity therefore gives
\[
\mathcal G\left(u(x_q+h),\sum_{\ell=1}^Q\llbracket u^\ell(x_q)+G^\ell(x_q)[h]\rrbracket\right)\leq C_Q[G]_{C^{0,\alpha}(U)}|h|^{1+\alpha}.
\]
By continuity of the paired jet, we may pass to the limit as $x_q\to x_0$, obtaining the same estimate with $x_q$ replaced by $x_0$.
This identifies the continuous paired jet as the Taylor first jet at every point.
Together with the H\"older continuity of $G$, this proves $u\in C^{1,\alpha}_{\mathrm{loc}}(B_1^n)$ in the sense of \cite[Definition~2.1]{Luc26}.

We now obtain the full norm estimate on $B_{1/2}^n$.
For every $x,y\in B_{1/2}^n$, the triangle inequality gives
\[
|Du(x)|\leq |Du(y)|+\mathcal G(Du(x),Du(y))\leq |Du(y)|+[Du]_{\alpha;B_{1/2}^n},
\]
since $\operatorname{diam}B_{1/2}^n=1$.
Integrating with respect to $y\in B_{1/2}^n$, dividing by $|B_{1/2}^n|$ and applying Cauchy--Schwarz, we obtain
\begin{equation}\label{eq:small-trace-gradient-sup}
\sup_{B_{1/2}^n}|Du|\leq C\|Du\|_{L^2(B_{1/2}^n)}+[Du]_{\alpha;B_{1/2}^n}.
\end{equation}
Similarly, integration along the segment joining $x$ and $y$ gives
\[
|u(x)|\leq |u(y)|+\mathcal G(u(x),u(y))\leq |u(y)|+\sup_{B_{1/2}^n}|Du|\,|x-y|.
\]
Integrating again with respect to $y$ yields
\begin{equation}\label{eq:small-trace-value-sup}
\sup_{B_{1/2}^n}|u|\leq C\|u\|_{L^2(B_{1/2}^n)}+\sup_{B_{1/2}^n}|Du|.
\end{equation}
Combining \eqref{eq:small-trace-gradient-sup} and \eqref{eq:small-trace-value-sup} with \eqref{eq:small-trace-holder} proves \eqref{eq:small-trace-full-norm}.
\end{proof}

\section{Applications to stationary maps}
\label{sec:area-core}

\begin{corollary}[$W^{2,2}$ Dirichlet-stationary maps]
\label{cor:W22-regularity}
For every $n\ge2$ and $k,Q\ge1$, a Dirichlet-stationary map
$u\in W^{2,2}(B^n_1,\mathcal A_Q(\R^k))$
has a $C^{1,\alpha}_{\mathrm{loc}}$ representative, with the exponent
and estimates of Theorem~\ref{thm:small-trace-regularity}.
\end{corollary}

\begin{proof}
By \cite[Lemma~3.10]{Luc26}, $\Delta u^\ell=0$ a.e.
Theorem~\ref{thm:small-trace-regularity} applies directly.
\end{proof}

For the area functional, write, for $P\in \R^{k\times n}$,
\[
 g(P):=I_n+P^TP,\qquad J(P):=\sqrt{\det g(P)},\qquad
 A(P):= D_PJ(P)=J(P)P g(P)^{-1}.
\]
The entries of $g(P)^{-1}$ are denoted by $g^{ij}(P)$, and we write
$A^a_i(P)=\frac{\partial J}{\partial P_i^a}(P)$ for the components
of the area flux.

\begin{lemma}[Non-divergence form of the $Q$-valued minimal surface system]
\label{lem:area-strong-equation}
Let $u\in W^{2,2}_{\mathrm{loc}}(\Omega,\mathcal A_Q(\R^k))\cap C^{0,1}_{\mathrm{loc}}(\Omega, \mathcal A_Q(\R^k))$ have a stationary graph varifold. Then, almost everywhere,
\begin{equation}
\label{eq:area-strong-equation}
 \sum_{i,j=1}^n g^{ij}(Du^\ell)D_{ij}u^\ell=0\in\R^k
 \qquad\text{for every }\ell.
\end{equation}
In particular,
\begin{equation}
\label{eq:area-small-trace}
 |\Delta u^\ell|\le|Du^\ell|^2|D^2u^\ell|,
 \qquad |\Delta u|\le\|Du\|_{L^\infty}^2|D^2u|.
\end{equation}
The last inequality may be applied on any region where $Du$ is bounded.
\end{lemma}

\begin{proof}
Let $U\Subset\Omega$. The area of the graph over $U$ is
\[
 \mathcal A(u;U)=\int_U\sum_{\ell=1}^Q J(Du^\ell).
\]
For $\varphi\in C_c^\infty(U)$ and
$\zeta\in C_c^\infty(\R^k;\R^k)$, the vertical ambient field
$X(x,y)=(0,\varphi(x)\zeta(y))$ induces the variation
$u_t=\sum_\ell\llbracket u^\ell+t\varphi\,\zeta(u^\ell)\rrbracket$.
Stationarity therefore gives
\begin{equation}
\label{eq:area-first-variation}
\begin{aligned}
 0
 &=\left.\frac{d}{dt}\right|_{t=0}\mathcal A(u_t;U)\\
 &=\int_U\sum_{\ell,i,a}
 \frac{\partial J}{\partial P_i^a}(Du^\ell)
 D_i\!\left(\varphi\,\zeta^a(u^\ell)\right)\\
 &=\int_U\sum_{\ell,i,a}
 A^a_i(Du^\ell)
 D_i\!\left(\varphi\,\zeta^a(u^\ell)\right)\\
 &=-\int_U\varphi\sum_{\ell,a}
 \bigl(\operatorname{div}[A(Du)]\bigr)^{\ell,a}\,\zeta^a(u^\ell).
\end{aligned}
\end{equation}

Define the $Q$-valued field
$f:=\operatorname{div} [A(Du)]=\sum_\ell\llbracket f^\ell\rrbracket$,
with the divergence understood through the approximate differential:
\[
 f^{\ell,a}:=\sum_{i,j,b}
 \frac{\partial A^a_i}{\partial P_j^b}(Du^\ell)
 D_{ij}u^{\ell,b},
\]
where
$j_2u=\sum_\ell\llbracket(u^\ell,Du^\ell,D^2u^\ell)\rrbracket$
is the approximate second jet of $u$.
Since $Du$ is locally bounded,
$f\in L^2_{\mathrm{loc}}(\Omega,\mathcal A_Q(\R^k))$.
The last equality in \eqref{eq:area-first-variation} follows from
the chain rule for the paired jet and integration by parts in
the symmetric sum.

Since $\varphi$ is arbitrary, a countable dense family of $\zeta$,
followed by separation of the distinct values of $u(x)$, shows that,
on one common set of full measure,
$\sum_{\{\ell:u^\ell(x)=z\}}f^\ell(x)=0$ for every
$z\in\operatorname{spt}u(x)$.
By Sobolev locality for $u$ and then for $j_1u$, coincident values
have equal first and second differentials almost everywhere, as in
\cite[Lemma~3.10]{Luc26}. Thus all $f^\ell$ within a coincidence
group agree, and each vanishes. Consequently,
$f=Q\llbracket0\rrbracket$ almost everywhere.

Fix a point $x$ where the approximate second jet exists and $f^\ell(x)=0$ for every $\ell$. For each $\ell$, choose a quadratic polynomial $w:\R^n\to\R^k$ whose second-order jet at $x$ is $(u^\ell(x),Du^\ell(x),D^2u^\ell(x))$. Set $\Sigma:=\operatorname{graph}w$ and $p:=(x,w(x))$, and denote its mean curvature vector by $\mathbf H$. The first variation formula (cf.~\cite[Chapter~2, \S\S~4--5]{Sim18}) gives
\[
J(Dw(x))\mathbf H^{\mathrm{vert}}(p)=\operatorname{div} [A(Dw)](x)=f^\ell(x)=0,
\]
where $\mathbf H^{\mathrm{vert}}(p)\in\R^k$ denotes the last $k$ components of $\mathbf H(p)\in\R^{n+k}$. Thus $\mathbf H(p)$ is horizontal. Being also normal to the graph, it must be zero. Moreover, expressing mean curvature as the trace of the second fundamental form gives
\[
0=\mathbf H(p)=\operatorname{proj}_{(T_p\Sigma)^\perp}\left(0,\sum_{i,j}g^{ij}(Dw(x))D_{ij}w(x)\right).
\]
The vector inside the projection is therefore both vertical and tangent to the graph, hence zero. Since the second-order jets agree, we obtain
\[
\sum_{i,j}g^{ij}(Du^\ell(x))D_{ij}u^\ell(x)=0.
\]
This holds almost everywhere for every $\ell$, proving \eqref{eq:area-strong-equation}.
Equivalently,
\[
 \Delta u^\ell
 =\sum_{i,j}a^{ij}(Du^\ell)D_{ij}u^\ell,
 \qquad a^{ij}(P):=\delta^{ij}-g^{ij}(P).
\]
Since $g(P)=I_n+P^TP$,
\[
 a(P)
 =I_n-g(P)^{-1}
 =\bigl(g(P)-I_n\bigr)g(P)^{-1}
 =P^TP\,g(P)^{-1}.
\]
Thus, using $\|g(P)^{-1}\|_{\mathrm{op}}\le1$, we have
$|a(P)|\le|P^TP|\le|P|^2$.
The asserted estimates follow by Cauchy--Schwarz and summation
over the sheets.
\end{proof}

\begin{corollary}[Small-slope area-stationary regularity]
\label{thm:area-regularity}
Let $u\in W^{2,2}(B_1^n,\mathcal A_Q(\R^k))\cap C^{0,1}_{\mathrm{loc}}(B_1^n,\mathcal A_Q(\R^k))$ have a stationary graph varifold.
If $\|Du\|_{L^\infty(B_1^n)}\le\sqrt{\delta_*(n,k,Q)}$, with $\delta_*$ as in Theorem~\ref{thm:small-trace-regularity}, then $u$ admits a $C^{1,\alpha}_{\mathrm{loc}}$ representative, with the exponent and interior estimates of that theorem.
In particular,
\begin{equation}
\label{eq:area-holder}
[Du]_{C^{0,\alpha}(B_{1/2}^n)}\le C\|Du\|_{L^\infty(B_1^n)},
\end{equation}
where $C=C(n,k,Q)<\infty$.
\end{corollary}

\begin{proof}
Lemma~\ref{lem:area-strong-equation} gives the small-Laplacian hypothesis of Theorem~\ref{thm:small-trace-regularity}, which therefore applies locally after translation and rescaling.

To prove \eqref{eq:area-holder}, we first establish a Caccioppoli estimate.
For $\eta\in C_c^1(B_1^n)$, testing \eqref{eq:weak-hodge} with $\Psi(x,Y)=\eta(x)^2Y$ and using the small-Laplacian bound gives
\[
(1-\delta_*^2)\int_{B_1^n}\eta^2|D^2u|^2\le 2(1+\delta_*)\int_{B_1^n}|\eta|\,|D\eta|\,|Du|\,|D^2u|.
\]
Since $\delta_*\le1/2$, Young's inequality yields
\begin{equation}
\label{eq:area-caccioppoli}
\int_{B_1^n}\eta^2|D^2u|^2\le C\int_{B_1^n}|D\eta|^2|Du|^2.
\end{equation}
Choosing a cutoff equal to one on $B_{3/4}^n$ and supported in $B_1^n$, with $|D\eta|\le C$, gives
\[
\|D^2u\|_{L^2(B_{3/4}^n)}\le C\|Du\|_{L^2(B_1^n)}\le C\|Du\|_{L^\infty(B_1^n)}.
\]
For each $x\in B_{1/2}^n$, apply \eqref{eq:small-trace-holder} to the rescaled map $z\mapsto4u(x+z/4)$, noting that $B_{1/4}^n(x)\subset B_{3/4}^n$.
Together with the preceding bound, this gives the required H\"older estimate for pairs of points in $B_{1/2}^n$ at distance less than $1/8$.
For pairs at distance at least $1/8$, the same estimate follows directly from the $L^\infty$ bound for $Du$.
This proves \eqref{eq:area-holder}.
\end{proof}

The Sobolev assumption in this corollary can also be obtained from a
geometric hypothesis. Since $DA(0)=I_{kn}$ and $DA$ is continuous and
self-adjoint, choose $\varepsilon_0=\varepsilon_0(n,k)\in(0,1/2]$ such that
\begin{equation}
\label{eq:area-ellipticity}
 \tfrac12|\dot P|^2
 \le DA(P)[\dot P]:\dot P
 \le2|\dot P|^2,
 \qquad |DA(P)[\dot P]|\le2|\dot P|
 \quad\text{if }|P|\le\varepsilon_0,\ \dot P\in\R^{k\times n}.
\end{equation}
Decrease $\delta_*$, if necessary, so that $\delta_*\le\varepsilon_0^2$;
all preceding conclusions remain valid.

\begin{lemma}[Capacitary extension of the first jet]
\label{lem:area-capacity}
Let $u:\Omega\to\mathcal A_Q(\R^k)$ be locally Lipschitz with
$\|Du\|_{L^\infty(\Omega)}\le\varepsilon_0(n,k)$.
Suppose that there is a relatively closed set $Z\subset\Omega$
with vanishing Sobolev $2$-capacity,
\begin{equation}
\label{eq:area-singular-set}
 \operatorname{cap}_2(Z)=0,
\end{equation}
and that $u$ admits local smooth minimal selections on
$\Omega\setminus Z$.
Then $u\in W^{2,2}_{\mathrm{loc}}(\Omega,\mathcal A_Q(\R^k))$.
Moreover, for every $\eta\in C_c^1(\Omega)$,
\begin{equation}
\label{eq:area-capacity-caccioppoli}
 \int_\Omega\eta^2|D^2u|^2
 \le C(n,k)\int_\Omega|D\eta|^2|Du|^2.
\end{equation}
\end{lemma}

\begin{proof}
Let $K\subset Z$ be compact and let $U\Subset\Omega$ be an
open neighbourhood of $K$.
By \eqref{eq:area-singular-set}, after passing to a
subsequence, we may choose
\begin{equation}
\label{eq:area-capacity-cutoffs}
 \zeta_m\in C_c^\infty(U),\qquad
 0\le\zeta_m\le1,\qquad
 \zeta_m=1\text{ near }K,\qquad
 \|\zeta_m\|_{W^{1,2}(U)}\longrightarrow0,
 \quad \zeta_m\longrightarrow0\text{ a.e.}
\end{equation}

Apply this construction to $K=Z\cap\operatorname{spt}\eta$.
For now, second derivatives are considered only on $\Omega\setminus Z$.
On each smooth selection neighborhood the minimal system gives
$\operatorname{div} [A(Du^\ell)]=0$. Differentiating in the $r$th
spatial direction yields the $\R^k$-valued system
\[
 \operatorname{div}\!\left(
 DA(Du^\ell)[D(D_ru^\ell)]
 \right)=0.
\]
Test with the $\R^k$-valued field
$\eta^2(1-\zeta_m)^2D_ru^\ell$, and sum over $r$ and $\ell$.
These tests have compact support away from $Z$.
A finite partition of unity subordinate to smooth selection neighborhoods
justifies the calculation: the summed expressions agree on overlaps by
locality, so the partition derivatives cancel.
The full matrix coercivity and boundedness in
\eqref{eq:area-ellipticity}, together with Young's inequality, give
\[
 \int\eta^2(1-\zeta_m)^2|D^2u|^2
 \le C(n,k)\int(1-\zeta_m)^2|D\eta|^2|Du|^2
 +C(n,k)\int\eta^2|D\zeta_m|^2|Du|^2.
\]
The last term tends to zero by the local boundedness of $Du$ and
\eqref{eq:area-capacity-cutoffs}. Fatou's lemma proves
\eqref{eq:area-capacity-caccioppoli}, initially for the Hessians off $Z$.
In particular these Hessians are locally square integrable up to $Z$.

The approximate paired jet $j_1u$ is locally bounded and, off $Z$, belongs
to $W^{1,2}_{\mathrm{loc}}$ with
$|D(j_1u)|^2=|Du|^2+|D^2u|^2$.
Let $\xi$ be an Almgren embedding and set $v:=\xi\circ j_1u$.
Testing its derivative identity off $Z$ with $(1-\zeta_m)\psi$,
where $\psi\in C_c^\infty(\Omega)$ and the cutoffs are chosen for
$K=Z\cap\operatorname{spt}\psi$, gives
\[
 \int_{\Omega\setminus Z}(1-\zeta_m)v\,D_i\psi
 =-\int_{\Omega\setminus Z}(1-\zeta_m)\psi\,D_iv
  +\int_{\Omega\setminus Z}\psi\,v\,D_i\zeta_m.
\]
The last term is bounded in absolute value by
\[
 \|\psi\|_\infty
 \|v\|_{L^2(\operatorname{spt}\psi)}
 \|D\zeta_m\|_{L^2}\longrightarrow0.
\]
The other terms converge by dominated convergence.
Since $Z$ has Lebesgue measure zero, this proves that
$v\in W^{1,2}_{\mathrm{loc}}(\Omega)$.
Hence $j_1u\in W^{1,2}_{\mathrm{loc}}$ and $u\in W^{2,2}_{\mathrm{loc}}$
by \cite[Definition~2.8]{Luc26}. Locality identifies its second differential
with the smooth Hessians off $Z$.
\end{proof}

We conclude by combining the preceding regularity results with
Bellettini's Lipschitz sheeting theorem for singular stable immersions
\cite[Theorems~4 and~5]{Bel25}. This gives sheeting with
$C^{1,\alpha}$ regularity in the multiple-valued sense, allowing branching
and crossings.

We first need the following definition, which extends the usual
notion of vanishing Sobolev $2$-capacity to the setting of
immersed manifolds.

\begin{definition}[$2$-capacity relative to a rectifiable varifold]
\label{def:relative-capacity}
Let $V$ be a rectifiable $n$-varifold in an open set $U\subset\R^{n+1}$. For $\chi\in C^1(U)$, write $\nabla^V\chi(x):=\operatorname{proj}_{T_xV}D\chi(x)$, where $T_xV$ is the approximate tangent plane, defined for $\|V\|$-almost every $x$. For compact $K\Subset U$, set
\[
c_{2,V}(K):=\inf_{\substack{\chi\in C^1(U),\ 0\le\chi\le1\\ \chi=1\text{ near }K}}\int_U\left(\chi^2+|\nabla^V\chi|^2\right)\,d\|V\|.
\]
For an arbitrary set $E\subset U$, define
\[
\operatorname{cap}_{2,V}(E):=\inf_{\substack{E\subset O\subset U\\ O\text{ open}}}\ \sup_{\substack{K\subset O\\ K\text{ compact}}}c_{2,V}(K).
\]
\end{definition}

This capacity is countably subadditive (see \cite[Section~4.7.1]{EG15} for the analogous argument in the Euclidean case). For relatively closed $Z\subset U$, the condition $\operatorname{cap}_{2,V}(Z)=0$ is equivalent to the existence of $\chi_j\in C^1(U)$, with $0\le\chi_j\le1$ and $\chi_j=1$ in a neighbourhood of $Z$, such that
\begin{equation}
\label{eq:vanishing-capacity}
\int_U\left(\chi_j^2+|\nabla^V\chi_j|^2\right)\,d\|V\|\longrightarrow0.
\end{equation}
This condition clearly implies $\|V\|(Z)=0$.
If $V_M:=\iota_\#|M|$ is induced by a smooth immersion
$\iota:M\to U\setminus Z$ and extended by zero across $Z$,
the area formula gives
\[
 \int_U\left(\chi^2+|\nabla^{V_M}\chi|^2\right)\,d\|V_M\|
 =
 \int_M\left(
 (\chi\circ\iota)^2+
 |\nabla^g(\chi\circ\iota)|_g^2
 \right)\,d\mathcal H_M^n
 \qquad\text{for every }\chi\in C^1(U),
\]
where $g=\iota^*g_{\mathrm{eucl}}$ and
$\mathcal H_M^n$ is the induced volume measure.

\medskip

For the next statement, we set the following notation: $C_R^{n+1}:=B_R^n(0)\times\mathbb R$, $\omega_n:=\mathcal L^n(B_1^n(0))$ and $E_M(R)^2:=R^{-n}\int_{C^{n+1}_R} \bigl(1-(\nu_M\cdot e_{n+1})^2\bigr)\,d\|V_M\|$, 
where $\nu_M$ is a unit normal to the immersion.

\begin{corollary}[$C^{1,\alpha}$ sheeting for singular stable immersions]
\label{cor:stable-sheeting}
Fix $\Lambda\ge1$.
There exist $\varepsilon=\varepsilon(n,\Lambda)>0$, $\alpha=\alpha(n,\Lambda)\in(0,1)$ and $C=C(n,\Lambda)<\infty$ with the following property.
Let $M$ be a smooth, properly immersed, two-sided stable minimal hypersurface without boundary in $C^{n+1}_2\setminus\Sigma$, where $\Sigma$ is relatively closed in $C^{n+1}_2$.
Assume that $\operatorname{cap}_{2,V_M}(\Sigma)=0$, $0\in\overline M$, $\|V_M\|(C^{n+1}_1)\le\Lambda$, $\sup_{M\cap C^{n+1}_2}|x_{n+1}|<\frac12$, and $E_M(1)\le\varepsilon$.
Then there are an integer $Q\ge1$ and a map $u\in C^{1,\alpha}(B_{1/4}^n,\mathcal A_Q(\R)) \cap W^{2,2}(B_{1/4}^n,\mathcal A_Q(\R))$ such that
\begin{equation}\label{eq:sheeting-graph}
\mathbf v(u)=V_M\mres C^{n+1}_{1/4}
\end{equation}
and
\begin{equation}\label{eq:sheeting-estimate}
\|Du\|_{C^{0,\alpha}(B_{1/4}^n)}\le CE_M(1).
\end{equation}
\end{corollary}

\begin{proof}
Let $\chi_j$ be cutoffs as in \eqref{eq:vanishing-capacity}.
For every $\Phi\in C_c^1(C^{n+1}_2;\R^{n+1})$, properness implies that the pullback of $(1-\chi_j)\Phi$ has compact support in the domain of the immersion.
Minimality and Cauchy--Schwarz therefore give
\[
\delta V_M(\Phi)=\lim_{j\to\infty}\delta V_M((1-\chi_j)\Phi)=0.
\]
Thus $V_M$ is stationary across $\Sigma$.

Choose $\varepsilon$ below the dimensional threshold in \cite[Theorems~4 and~5]{Bel25}, which apply under the present capacitary assumption (see the discussion at the beginning of Part~II of \cite{Bel25}; see also \cite[Remark~1.5]{WZ26}).
Those results give a Lipschitz $Q$-valued map $\widetilde u$ on $B_{1/2}^n$ such that
\[
\mathbf v(\widetilde u)=V_M\mres C_{1/2}^{n+1}.
\]
Each ordered scalar selection has Lipschitz constant at most $C(n)E_M(1)$.
The area formula and the mass bound give
\[
Q\omega_n2^{-n}\le\|V_M\|(C^{n+1}_{1/2})\le\Lambda,
\]
so only finitely many values of $Q$ can occur, in terms of $n,\Lambda$.
In the summed norm used throughout this paper,
\[
\|D\widetilde u\|_{L^\infty(B_{1/2}^n)}
\le\sqrt Q\,C(n)E_M(1).
\]
Decrease $\varepsilon(n,\Lambda)$ so that this is at most the small-slope threshold $\sqrt{\delta_*(n,1,Q)}$ from Corollary~\ref{thm:area-regularity}, for every possible $Q$.

Let $\pi:\R^{n+1}\to\R^n$ be the horizontal projection and set $Z:=\pi(\Sigma\cap\operatorname{graph}\widetilde u)$.
The continuity of $\widetilde u$ shows that $Z$ is relatively closed in $B_{1/2}^n$.
Define
\[
\zeta_j(x):=\max_{y\in\operatorname{spt}\widetilde u(x)}\chi_j(x,y).
\]
These functions are locally Lipschitz, satisfy $0\le\zeta_j\le1$, and equal $1$ in a neighbourhood of $Z$.
The chain rule, the Lipschitz bounds for the ordered selections, and the area formula give
\[
\|\zeta_j\|_{W^{1,2}(B_{1/2}^n)}\longrightarrow0.
\]
Localization and smoothing therefore yield $\operatorname{cap}_2(Z)=0$.

By the local graphical description in \cite[Theorem~5]{Bel25}, $\widetilde u$ admits local smooth minimal selections on $B_{1/2}^n\setminus Z$.
Since $\sqrt{\delta_*(n,1,Q)}\le\varepsilon_0(n,1)$, Lemma~\ref{lem:area-capacity} gives $\widetilde u\in W^{2,2}_{\mathrm{loc}}(B_{1/2}^n,\mathcal A_Q(\R))$.
Its graph varifold is stationary, so Corollary~\ref{thm:area-regularity}, applied after translation and rescaling on interior balls, gives $\widetilde u\in C^{1,\alpha}_{\mathrm{loc}}(B_{1/2}^n,\mathcal A_Q(\R))$.
Taking the minimum of the exponents, and the maximum of the estimate constants, over the finitely many possible values of $Q$ makes them depend only on $n,\Lambda$.

Applying \eqref{eq:area-holder} after rescaling on $B_\rho^n$, with $\rho<1/2$, and letting $\rho\uparrow1/2$ gives
\[
[D\widetilde u]_{C^{0,\alpha}(B_{1/4}^n)}
\le C\|D\widetilde u\|_{L^\infty(B_{1/2}^n)}
\le CE_M(1).
\]
Set $u:=\widetilde u|_{B_{1/4}^n}$.
Since $B_{1/4}^n\Subset B_{1/2}^n$, the local regularity established above gives
\[
u\in C^{1,\alpha}(B_{1/4}^n,\mathcal A_Q(\R))
\cap W^{2,2}(B_{1/4}^n,\mathcal A_Q(\R)).
\]
The graphical representation gives \eqref{eq:sheeting-graph}, and the preceding estimates give \eqref{eq:sheeting-estimate}.
\end{proof}

\section{Regularity and conditional compactness for stable minimal immersions}
\label{sec:stable-immersions}

In this section, we show how Corollary~\ref{cor:stable-sheeting},
combined with \cite[Theorem~5.7]{Luc26}, complements the regularity
and precompactness theory of Wang and Zhang \cite{WZ26}.
Our purpose is to explain the resulting control of the singular set
of limits and to identify the remaining obstruction to compactness.

\begin{definition}[The class $\mathcal V(U)$]
\label{def:stable-immersion-class}
Let $n\ge2$ and let $U\subset\R^{n+1}$ be open.
We denote by $\mathcal V(U)$ the class of integral $n$-varifolds $V$ in $U$ for which there exist a relatively closed set $Z\subset U$, a possibly disconnected smooth $n$-manifold $M$ without boundary, and a smooth, proper, two-sided immersion
\[
\iota:M\longrightarrow U\setminus Z.
\]
We equip $M$ with the induced metric $g=\iota^*g_{\mathrm{eucl}}$ and require the following properties.
\begin{enumerate}
\item
The varifold $V$ is induced by the immersion: $V=\iota_\#|M|$,
extended by zero across $Z$, where $|M|$ denotes the multiplicity-one varifold associated with $(M,g)$.

\item
The immersion is stationary: for every smooth variation $(\iota_t)$ with $\iota_0=\iota$ and $\iota_t\equiv\iota$ outside a fixed compact set $K\subset M$, we have
\[
\left.\frac{d}{dt}\right|_{t=0}
\bigl\|(\iota_t)_\#|M|\bigr\|(W)=0,
\]
where $W\Subset U\setminus Z$ is any open set containing $\iota_t(K)$ for all sufficiently small $t$.
Equivalently, the mean curvature vector of $\iota$ vanishes.

\item
The immersion is stable: for every such variation and choice of $W$, we have
\[
\left.\frac{d^2}{dt^2}\right|_{t=0}
\bigl\|(\iota_t)_\#|M|\bigr\|(W)\geq0.
\]
Equivalently,
\[
\int_M |A_\iota|_g^2\varphi^2\,d\mathcal H_M^n
\leq
\int_M |\nabla^g\varphi|_g^2\,d\mathcal H_M^n
\qquad\text{for every }\varphi\in C_c^1(M),
\]
where $A_\iota$ is the second fundamental form and $\mathcal H_M^n$ is the measure associated with the induced metric.

\item
$Z$ has vanishing $2$-capacity relative to $V$ in the sense
of Definition~\ref{def:relative-capacity}:
$\operatorname{cap}_{2,V}(Z)=0$.
\end{enumerate}
\end{definition}

Note that multiplicities are included by allowing repeated components of $M$, and the zero varifold is included by allowing $M=\varnothing$.
Note also that, compared with the formulation in \cite{Bel25,WZ26}, we replace local finiteness of $\mathcal H^{n-2}(Z)$ with vanishing $2$-capacity relative to $V$.
The results of Wang and Zhang remain applicable under this weaker hypothesis \cite[Remark~1.5]{WZ26}.

Let $\overline{\mathcal V}(U)$ denote the class of varifold limits of sequences in $\mathcal V(U)$ with locally uniformly bounded mass.
The compactness problem is to prove that
\[
\overline{\mathcal V}(U)=\mathcal V(U).
\]
We use the immersed notion of regularity.
A point $x\in\operatorname{spt}\|V\|$ is regular if, in some neighbourhood, $V$ is a finite integer sum of relatively closed smooth\footnote{Here smooth means $C^\infty$, not merely $C^1$; this technical choice allows us to establish the minimality of the individual sheets directly from stationarity in Step~1 of the proof of Theorem~\ref{thm:conditional-compactness}.} embedded hypersurfaces without interior boundary.
Crossings of smooth sheets are allowed.
We write $\operatorname{sing}V$ for the complement of this regular set in $\operatorname{spt}\|V\|$.

Combining the sheeting theory of Wang and Zhang \cite{WZ26} with Corollary~\ref{cor:stable-sheeting} and the branch-set dimension bound of \cite[Theorem~5.7]{Luc26}, we obtain
\[
\dim_{\mathcal H}\operatorname{sing}V\le n-2
\qquad\text{for every }V\in\overline{\mathcal V}(U).
\]

This estimate, however, does not by itself imply
vanishing $2$-capacity and therefore does not yet establish
compactness of $\mathcal V(U)$.
The remaining obstruction can be isolated in the following
conjecture, which records a natural expectation in the
regularity theory of branched minimal immersions.

\medskip
\noindent
\textbf{Conjecture $(\star)$ (zero capacity of branch sets).}
For every $n\ge2$, $Q\ge1$, $\alpha\in(0,1)$ and open set $\Omega\subset\R^n$, let $u\in C^{1,\alpha}_{\mathrm{loc}}(\Omega;\mathcal A_Q(\R))$ have stationary graph varifold.
Then its branch set $\mathcal B_u$\footnote{Here $\mathcal B_u$ consists of the points at which $u$ does not admit a local representation by $Q$ $C^1$ single-valued maps.} has vanishing Sobolev $2$-capacity in $\Omega$.

\medskip
The conjecture already holds in the $2$-valued case by the work of Krummel and Wickramasekera \cite[Theorem~C]{KW21}, who prove that the branch set of a stationary $C^{1,\alpha}$ $2$-valued graph is countably $(n-2)$-rectifiable.
This implies $\sigma$-finiteness of its $\mathcal H^{n-2}$-measure and hence vanishing $2$-capacity.

Assuming the conjecture, Corollary~\ref{cor:stable-sheeting}, together with the results of Wang and Zhang,
yields sequential compactness of $\mathcal V(U)$
under locally uniform mass bounds. We prove this conditional result in Theorem~\ref{thm:conditional-compactness} below.

We first recall the following classical lemma, which we prove for completeness.

\begin{lemma}[Tangential contacts]
\label{lem:tangential-contact-capacity}
Let $\Omega\subset\R^n$ be connected and open, and let
$u,v\in C^\infty(\Omega)$ be distinct solutions of the
minimal surface equation. Then
\[
 C:=\{x\in\Omega:u(x)=v(x),\ Du(x)=Dv(x)\}
\]
is contained in the union of countably many smooth submanifolds of dimension $n-2$.
In particular, $C$ is countably $(n-2)$-rectifiable and has locally vanishing Sobolev $2$-capacity.
\end{lemma}

\begin{proof}
We follow the argument of \cite[Lemma~1.9]{HS89}, which we recall here for completeness.

By the analyticity of smooth solutions of the minimal
surface equation \cite{Mor58}, the function $w=u-v$
is real analytic. Moreover, it satisfies a PDE of the form
\[
\sum_{ij} a^{ij}D_{ij}w+ \sum_i b^iD_iw=0,
\]
with smooth coefficients and locally uniformly elliptic
principal part. Since $w\not\equiv0$, every point of $C$
has finite vanishing order. Write
\[
 C=\bigcup_{r\ge2}C_r,
 \qquad
 C_r:=\{p:D^\gamma w(p)=0\text{ for }|\gamma|<r,\
                         D^rw(p)\ne0\}.
\]
Let $P_r$ be the nonzero leading homogeneous Taylor
polynomial at $p\in C_r$. As $w(p+tz)/t^r\to P_r$ locally smoothly, $\sum_{ij} a^{ij}(p)D_{ij}P_r=0$.
The vectors $\{D (D^\beta w)(p):|\beta|=r-1\}$ span a space of dimension at least two. Indeed, otherwise,
after rotating coordinates, they would all be parallel to $e_1$.
Consequently, $D_jD^\beta w(p)=0$ for every $j\ge2, |\beta|=r-1$.
Since derivatives commute, every derivative of order $r$
involving a coordinate other than the first would vanish.
Thus $P_r(z)=c z_1^r$ for some $c\ne0$ and $r\geq 2$, contradicting ellipticity.

Choose two independent such vectors. The implicit function
theorem shows that the common zero set of the corresponding
derivatives is locally a smooth submanifold of codimension
two containing $C_r$ in a neighbourhood of $p$.
A countable subcover for each $r$ proves that $C$ is contained in the union of countably many smooth manifolds of dimension $n-2$.
When $n\geq 3$, the capacity conclusion follows from $\sigma$-finiteness
of $\mathcal H^{n-2}$, \cite[Theorem~4.16]{EG15}, and
countable subadditivity. When $n=2$, use countability and
the zero $2$-capacity of a point (which can be checked using the logarithmic cutoff trick).
\end{proof}

\begin{theorem}[Conditional compactness]
\label{thm:conditional-compactness}
Let $n\ge2$, let $U\subset\R^{n+1}$ be open, and let $(V_j)_{j\geq 1}$ be a sequence in $\mathcal V(U)$ with uniformly locally bounded masses. Then, after passing to a subsequence, $V_j$ converges as varifolds to a stationary integral $n$-varifold $V$ in $U$, and there exist a smooth $n$-manifold $M$ without boundary and a smooth, proper, two-sided stable minimal immersion
\[
 \iota:M\longrightarrow U\setminus\operatorname{sing}V
\]
such that
\[
 \iota_\#|M|
 =V\mathbin{\llcorner}
       (U\setminus\operatorname{sing}V).
\]
The singular set of the limit satisfies
\[
 \dim_{\mathcal H}\operatorname{sing}V\le n-2.
\]
Moreover, if Conjecture~$(\star)$ holds, then $V\in \mathcal V(U)$\footnote{For $n=2$, the compactness conclusion is obtained
independently by Stanbury in \cite[Theorem~1.17]{Sta26},
without assuming Conjecture~$(\star)$. His proof uses a
calibration argument to establish discreteness of the
branch sets of the limiting graphical maps.}.
\end{theorem}

\begin{proof}
\medskip\noindent
\emph{Step 1: construction of the minimal immersion.}
We first extend stationarity across $Z_j$.
For each $j$, let $\iota_j:M_j\to U\setminus Z_j$ be an immersion representing $V_j$ as in Definition~\ref{def:stable-immersion-class}, so that $V_j=(\iota_j)_\#|M_j|$, extended by zero across $Z_j$.
Let $\chi_{j,r}$ be capacity cutoffs for $Z_j$.
Stationarity away from $Z_j$ and Cauchy--Schwarz give $\delta V_j(\Phi)  =\lim_{r\to\infty}\delta V_j((1-\chi_{j,r})\Phi)=0$ for every $\Phi\in C_c^1(U;\R^{n+1})$.
Thus each $V_j$ is stationary in $U$.
Allard compactness now gives a varifold-convergent subsequence with stationary integral limit $V$.

Recall that, by definition, every point of $\operatorname{reg}V$ has an ambient neighbourhood in which $V$ is a finite integer sum of relatively closed smooth embedded hypersurfaces without interior boundary.

The smooth local branches of $V$ on $\operatorname{reg}V$ are individually minimal. Indeed, write locally $V=\sum_{a=1}^s q_a|\Sigma_a|$.
For each $a$, the set
$\Sigma_a\setminus\bigcup_{b\ne a}
\partial_{\Sigma_a}(\Sigma_a\cap\Sigma_b)$
is open and dense in $\Sigma_a$. Near each of its
points, every other branch either coincides with
$\Sigma_a$ or is absent, so stationarity of $V$ gives
$\mathbf H_{\Sigma_a}=0$.
By continuity of $\mathbf H_{\Sigma_a}$, the same conclusion holds everywhere on $\Sigma_a$.

Let $\widehat M$ be the natural desingularization of $\operatorname{reg}V$, obtained by separating its distinct smooth minimal local branches.
More precisely, $\widehat M$ consists of pairs $(x,[\Sigma]_x)$, where $x\in\operatorname{reg}V$ and $[\Sigma]_x$ is the germ at $x$ of a smooth embedded minimal hypersurface $\Sigma$ passing through $x$, relatively closed in some ambient neighbourhood of $x$ and contained there in $\operatorname{spt}\|V\|$.

By unique continuation, any two local decompositions into distinct smooth minimal branches agree up to a permutation after restriction to a sufficiently small common neighbourhood.
Indeed, choose smooth functions $f,g_1,\ldots,g_s$ on a common connected base $B^n$ such that $\operatorname{graph}f$ represents an arbitrary minimal germ and $\operatorname{graph}g_1,\ldots,\operatorname{graph}g_s$ are the branches of a local decomposition.
Since $\operatorname{graph}f$ is contained in $\operatorname{spt}\|V\|$, we have $B^n=\bigcup_{b=1}^s\{f=g_b\}$.
At least one of these closed coincidence sets has nonempty interior, so unique continuation gives $f\equiv g_b$ on $B^n$ for some $b$.
This also shows that every minimal germ is represented by a branch of any local decomposition.

The local branches therefore provide compatible smooth charts on $\widehat M$.
Unique continuation also ensures that $\widehat M$ is Hausdorff, and a countable cover of $\operatorname{reg}V$ by neighbourhoods with finitely many branches gives a countable atlas on $\widehat M$, hence second countability.
Coincident germs are counted once, with their total multiplicity recorded in a locally constant positive integer function $q:\widehat M\to\mathbb N$.
The natural projection
\[
\widehat\iota:\widehat M\longrightarrow U\setminus\operatorname{sing}V,
\qquad (x,[\Sigma]_x)\longmapsto x,
\]
is a smooth, proper minimal immersion (properness follows from the local finiteness and relative closedness of the branches).
By construction, $\widehat\iota_\#(q|\widehat M|)=V\mres (U\setminus\operatorname{sing}V)$.
We equip $\widehat M$ with the induced metric $\widehat g:=\widehat\iota^*g_{\mathrm{eucl}}$.

\medskip\noindent
\emph{Step 2: $C^{1,\alpha}$ local graphical decompositions and convergence.}
We now combine the sheeting theorem of Wang and Zhang with our regularity estimates.
Consider a point at which $V$ admits a tangent cone $\sum_a q_a|P_a|$, where the $P_a$ are distinct hyperplanes whose common intersection has dimension at least $n-2$.
Choosing a sufficiently small fixed scale and then taking $j$ large, we may apply \cite[Theorem~1.9 and Remark~1.5]{WZ26} to the corresponding rescalings of $V_j$.
On a smaller neighbourhood, this gives a decomposition into stationary graph varifolds associated with Lipschitz $q_a$-valued maps $u_{j,a}$ over fixed balls $B_a^n\subset P_a$, with uniformly small slopes.

Choose the initial slope bounds sufficiently small that the tangent planes of the different graphical families lie in fixed neighbourhoods of the corresponding planes $P_a$ with pairwise disjoint closures.
Applying the argument in the proof of Corollary~\ref{cor:stable-sheeting} to each graphical component gives uniform interior $C^{1,\alpha}$ estimates.
By compactness, every subsequence admits a further subsequence along which all the maps $u_{j,a}$ converge locally in $C^1$.
Along any such subsequence, $C^1$ convergence implies convergence of the graph varifolds, so passing to the limit in $V_j=\sum_a\mathbf v(u_{j,a})$ shows that the limiting graph varifolds sum locally to $V$.
In a smaller fixed graphical neighbourhood, the graph varifold of the $a$-th limiting map is precisely the part of $V$ whose tangent planes lie in the chosen neighbourhood of $P_a$.
Thus each limiting map is uniquely determined by $V$, so all subsequential limits agree.
Consequently, $u_{j,a}\to u_a$ locally in $C^1$, without further extraction, for maps $u_a\in C^{1,\alpha}_{\mathrm{loc}}$.
Their graph varifolds are stationary.
In coordinates centred at the point where the reference tangent cone is taken, we have $j_1u_a(0)=q_a\llbracket(0,0)\rrbracket$.

Choose a unit vector $e$ transverse to all the reference hyperplanes $P_a$.
Since the limiting graphical families are tangent to the corresponding planes $P_a$ at the centre, shrinking the neighbourhood makes the Lipschitz constants of their ordered scalar selections as small as required.
By $C^1$ convergence, the same holds for $u_{j,a}$ for all sufficiently large $j$.
As each $\pi_{e^\perp}|_{P_a}$ is invertible, projection onto $e^\perp$ is therefore bi-Lipschitz on each of these scalar graphs, with constants independent of $j$.
Their images contain a common fixed base ball.
The resulting reparametrized $q_a$-valued maps are $C^{1,\alpha}$ with uniform bounds on smaller balls (see \cite[Lemma~5.4 and Lemma~5.5]{Luc26} and the arguments in their proofs).
They converge uniformly to the corresponding reparametrized limits, and hence in $C^1$ on a smaller fixed base ball by compactness.
We restrict the ambient neighbourhood to its intersection with the cylinder over this ball.
Thus, summing the reparametrized maps, we obtain maps $u_j,u\in C^{1,\alpha}(\widetilde B^n,\mathcal A_Q(\R))$ on a common fixed base ball $\widetilde B^n\subset e^\perp$, where $Q:=\sum_a q_a$.
These satisfy $V_j=\mathbf v(u_j)$ and $V=\mathbf v(u)$ in the chosen graphical neighbourhood for all sufficiently large $j$, and $u_j\to u$ in $C^1(\widetilde B^n)$.

Let $T$ be the complement in $\operatorname{spt}\|V\|$ of the union of these graphical neighbourhoods.
Then $T$ is relatively closed.
The cone analysis in \cite[Section~5]{WZ26} shows that every tangent cone of $V$ whose spine has dimension at least $n-2$ is a union of hyperplanes of the type considered above.
Here stability and two-sidedness are required only for the approximating immersions: by a diagonal argument, tangent cones of $V$ and their iterated tangent cones can be realised as limits of rescalings of the $V_j$.
Every point admitting such a tangent cone is therefore contained in one of the graphical neighbourhoods constructed above.
Consequently, for $n\ge3$, $T\subseteq\mathcal S_{n-3}(V)$, where $\mathcal S_m(V)$ denotes the standard stratum of points at which every tangent cone has spine dimension at most $m$.
Almgren's stratification theorem gives $\dim_{\mathcal H}T\le n-3$.
When $n=2$, every tangent cone has spine dimension at least $0=n-2$, so the preceding argument gives $T=\varnothing$.

Let $C_{\mathrm{tan}}$ be the set of tangential contacts between distinct smooth germs in the immersed regular part.

Let $x\in\operatorname{reg}V\setminus(T\cup C_{\mathrm{tan}})$.
We show that, in a neighbourhood of $x$, we can write $V=\sum_{b=1}^s q_b|\Sigma_b|$ and $V_j=\sum_{b=1}^s\mathbf v(u_j^b)$, where the $\Sigma_b$ are distinct smooth embedded minimal hypersurfaces through $x$ and the $q_b$-valued maps $u_j^b$ are defined on a common fixed base ball $B^n$, with $u_j^b\to u^b$ in $C^1(B^n)$ and $\mathbf v(u^b)=q_b|\Sigma_b|$.
Since $x\notin T$, the preceding construction gives an open neighbourhood $O$ of $x$ and maps $u_j,u\in C^{1,\alpha}(\widetilde B^n,\mathcal A_Q(\R))$ over a common fixed ball $\widetilde B^n$ in an $n$-dimensional plane $P$, such that $V_j=\mathbf v(u_j)$ and $V=\mathbf v(u)$ in $O$ for all sufficiently large $j$, and $u_j\to u$ in $C^1(\widetilde B^n)$.
Since $x\in\operatorname{reg}V$, after shrinking $O$ we have $V=\sum_{b=1}^s q_b|\Sigma_b|$ in $O$, where the $\Sigma_b$ are distinct smooth embedded minimal hypersurfaces through $x$.
Restrict $\widetilde B^n$ to a smaller ball $B^n$ centred at the projection of $x$, so that the sheets of $u$ passing through $x$ lie in $O$ and remain separated from the other sheets.
Discard these other sheets, together with the corresponding sheets of $u_j$ for all sufficiently large $j$.
Keeping the same notation, we obtain $u_j,u\in C^{1,\alpha}(B^n,\mathcal A_Q(\R))$, now with $Q:=\sum_{b=1}^s q_b$.
Then $V_j=\mathbf v(u_j)$ and $V=\mathbf v(u)$ near $x$, with $u_j\to u$ in $C^1(B^n)$ and $u=\sum_{b=1}^s u^b$, where $u^b:=q_b\llbracket f_b\rrbracket$ and $f_b$ is the smooth function parametrizing $\Sigma_b$ over $B^n$.
Since $x\notin C_{\mathrm{tan}}$, the differentials $Df_b$ are pairwise distinct at the projection of $x$ and remain separated after shrinking $B^n$.
The $C^1$ convergence therefore gives, for all sufficiently large $j$, a corresponding splitting $u_j=\sum_{b=1}^s u_j^b$, with $u_j^b\in C^{1,\alpha}(B^n,\mathcal A_{q_b}(\R))$ and $u_j^b\to u^b$ in $C^1(B^n)$ (see \cite[Corollary~2.7]{Luc26}).

Set
\[
\widehat E:=\widehat\iota^{-1}(T\cup C_{\mathrm{tan}}).
\]
The dimension bound for $T$, local graphical coordinates, and Lemma~\ref{lem:tangential-contact-capacity} show that $\widehat E$ has locally vanishing intrinsic $2$-capacity.
In particular, $\mathcal H^{n-1}(\widehat E)=0$.
Moreover, $\widehat E$ is closed in $\widehat M$, since $T$ is relatively closed in $\operatorname{spt}\|V\|$ and $C_{\mathrm{tan}}$ is relatively closed in $\operatorname{reg}V$ by the local finiteness, relative closedness and smoothness of the branches.

\medskip\noindent
\emph{Step 3: two-sidedness.}
We show that every one-sided connected component of $\widehat M$ has even multiplicity.

Fix $p\in\widehat M\setminus\widehat E$ and set $x:=\widehat\iota(p)\in \operatorname{reg}V \setminus (T\cup C_{\mathrm{tan}})$.
Use the decomposition obtained at the end of Step~2 in a sufficiently small neighbourhood of $x$, and label its branches so that $\Sigma_1$ is locally parametrized by $\widehat\iota$ near $p$.
Set $V_j^1:=\mathbf v(u_j^1)$ and $q:=q_1=q(p)$.
Thus $u_j^1\to u^1$ in $C^1(B^n,\mathcal A_q(\R))$, with $\mathbf v(u^1)=q|\Sigma_1|$.
Let $\pi_{B^n}$ denote the graphical projection.
Let $M_j^1\subset M_j$ be the open portion of the immersion domain corresponding to $u_j^1$, so that $(\iota_j)_\#|M_j^1|=V_j^1$.
Thus $\iota_j(M_j^1)=\operatorname{graph}u_j^1\setminus Z_j$.

We now show that, for all sufficiently large $j$, the numbers of positively and negatively oriented sheets of $M_j^1$ relative to a chosen continuous unit normal to $\Sigma_1$ are constant away from an exceptional set $S_j$ of vanishing $2$-capacity in $B^n$.
Set $S_j:=\pi_{B^n}(Z_j\cap\operatorname{graph}u_j^1)$ and $p_j:=\pi_{B^n}\circ\iota_j|_{M_j^1}$.
The set $S_j$ is relatively closed in $B^n$, and the restriction $p_j:M_j^1\setminus p_j^{-1}(S_j)\longrightarrow B^n\setminus S_j$ is a proper local diffeomorphism, hence a finite covering of degree $q$.
The cutoff construction in the proof of Corollary~\ref{cor:stable-sheeting}, applied to $u_j^1$, gives $\operatorname{cap}_2(S_j)=0$.
In particular, $B^n\setminus S_j$ is connected.
Fix a smooth global unit normal $\nu_j$ on each $M_j$, and choose a smooth unit normal $\nu$ to $\Sigma_1$, viewed as a function of the base point $y\in B^n$ ($\Sigma_1$ is orientable since it is a graph).
For all sufficiently large $j$, the $C^1$ convergence ensures that $\nu_j(z)\cdot\nu(y)\neq0$ whenever $y\in B^n\setminus S_j$ and $z\in p_j^{-1}(y)$.
Define $q_j^\pm(y):=\#\{z\in p_j^{-1}(y):\ \pm\nu_j(z)\cdot\nu(y)>0\}$.
On each local sheet of the covering, this scalar product is continuous and never vanishes, so its sign is locally constant.
The counts $q_j^\pm(y)$ are therefore locally constant.
Since $B^n\setminus S_j$ is connected, they are constant, and we denote their values by $q_j^\pm$.
Each fibre contains $q$ points, so $q_j^++q_j^-=q$.

We next compare the sign counts $q_{j,\alpha}^\pm,q_{j,\beta}^\pm$ of two graphical charts around $p\in\widehat M\setminus\widehat E$ whose reference normals $\nu_\alpha,\nu_\beta$ agree near $p$.
By the tangent-plane separation in Step~2 and $C^1$ convergence, in a sufficiently small fixed ambient neighbourhood $W$ of $\widehat\iota(p)$ both charts select the same portion $N_j\subset M_j$ for all sufficiently large $j$.
Set $v:=\nu_\alpha(p)=\nu_\beta(p)$ and $N_j^\pm:=\{z\in N_j:\ \pm\nu_j(z)\cdot v>0\}$; after shrinking $W$, these are the positive and negative portions identified in both charts.
The two charts thus describe the same immersed portions, away from the respective exceptional sets, by $q_{j,\alpha}^\pm$-valued and $q_{j,\beta}^\pm$-valued maps, but with respect to different graphical projections.
To show that $q_{j,\alpha}^\pm=q_{j,\beta}^\pm$ for all sufficiently large $j$, let $\Sigma\subset W$ be the common smooth limiting hypersurface.
In each chart, all sheets converge in $C^1$ to the single-valued graph representing $\Sigma$.
Choose a nonnegative $\phi\in C_c^\infty(W)$ with $\int_\Sigma\phi\,d\mathcal H^n=1$.
Since the exceptional sets have measure zero and the sign counts are constant and bounded by $q$, the area formula and $C^1$ convergence in the two charts give $\int_{N_j^\pm}\phi(\iota_j(z))\,d\mathcal H^n_{M_j}(z)=q_{j,\alpha}^\pm+o_j(1)=q_{j,\beta}^\pm+o_j(1)$.
Thus $q_{j,\alpha}^\pm-q_{j,\beta}^\pm$ is an integer tending to zero and therefore vanishes for all sufficiently large $j$.

Let $C$ be a one-sided connected component of $\widehat M$ (if there are none, we are done), and denote its constant multiplicity by $q$.
Its orientation double cover is connected and remains path connected after removing the lift of $\widehat E$, since $\mathcal H^{n-1}(\widehat E)=0$.
Thus $C\setminus\widehat E$ contains a loop along which a continuous unit normal reverses sign.
Subdivide this loop into finitely many segments contained in connected graphical charts in $C\setminus\widehat E$, using the same chart for the first and last segments.
Continue a unit normal along the loop and choose the reference normal on each segment accordingly.
By the preceding comparison, for all sufficiently large $j$ the constant sign counts agree at every transition, so the ordered pair remains unchanged along this chain of charts.
Denote the initial pair by $(q_j^+,q_j^-)$.
On returning to the initial chart, the reference normal has reversed, so the final pair is $(q_j^-,q_j^+)$.
Consequently, $q_j^+=q_j^-$ and $q=q_j^++q_j^-=2q_j^+$, so $q$ is even.

For each component of $\widehat M$ of multiplicity $q$, take $q$ copies of the component if it is two-sided, and $q/2$ copies of its orientation double cover otherwise.
The resulting covering $\pi:M\to\widehat M$ has degree $q$ over each component, and $\iota:=\widehat\iota\circ\pi:M\longrightarrow U\setminus\operatorname{sing}V$ is a smooth, proper, two-sided minimal immersion satisfying $\iota_\#|M|=V\mathbin{\llcorner}(U\setminus\operatorname{sing}V)$.

\medskip\noindent
\emph{Step 4: stability.}
We first prove the scalar stability inequality on $\widehat M$, which need not be two-sided.
Fix $\varphi\in C_c^\infty(\widehat M\setminus\widehat E)$ and extend it by zero to $\widehat M$.
We construct
\[
 Y\in C_c^\infty\bigl(
 (U\setminus \operatorname{sing}V)\times\mathbb S^n;
 \R^{n+1}\bigr)
\]
with spatial support disjoint from $T\cup C_{\mathrm{tan}}$,
odd in its second variable, and satisfying
\[
 Y(\widehat\iota(p),\nu^{\pm}(p))=\varphi(p)\nu^\pm (p)
\]
for both unit normals $\nu^\pm (p)$ at each $p\in\widehat M$.
We also require $D_\nu Y\equiv 0$ near every pair $(\widehat\iota(p),\nu^\pm (p))$ with $p\in\widehat M$.

To construct $Y$, cover $\widehat\iota(\operatorname{spt}\varphi)$ by sufficiently small ambient patches in $U\setminus \operatorname{sing} V$ and disjoint from $T\cup C_{\mathrm{tan}}$, on which distinct local branches of $V$ have separated tangent planes.
On each patch, let $\pi_a$ be the tubular projection onto the $a$-th branch, let $\varphi_a$ denote the function induced there by $\varphi$ on the corresponding lifted branch, and choose a smooth local unit normal $\nu_a$.
Choose $\theta_a\in C^\infty(\mathbb S^n)$ equal to $1$ near the values of $\nu_a$ along $\Sigma_a$ and supported in a spherical cap, with the caps chosen for the various indices $a$ and their antipodal images all pairwise disjoint.
Define
\[
Y_{\mathrm{loc}}(X,\nu)
:=\sum_a \varphi_a(\pi_a(X))\,\nu_a(\pi_a(X))
\bigl(\theta_a(\nu)-\theta_a(-\nu)\bigr).
\]
The separation of the caps gives the prescribed values at both unit normals to each branch, while the constant plateaus ensure that $D_\nu Y_{\mathrm{loc}}\equiv 0$ near all such pairs.
Combining these local fields with a spatial partition of unity subordinate to the patches and summing to $1$ near $\widehat\iota(\operatorname{spt}\varphi)$ gives the required compactly supported field $Y$, which is odd in $\nu$.

Choose global unit normals $\nu_j$ on $M_j$, extend $Y$ by zero to $U\times\mathbb S^n$, and set $Y_j:=Y(\iota_j,\nu_j)$.
Uniform convergence of the tangent planes implies that,
for all sufficiently large $j$,
\[
 D_\nu Y(\iota_j,\nu_j)=0,
 \qquad
 dY_j=D_XY(\iota_j,\nu_j)\,d\iota_j.
\]

To pass the stability inequality to the limit, we use the fact that, for linear variations of an immersion, the second variation integrand depends continuously on the tangent plane and quadratically on the differential of the variation field (see \cite[Chapter~2, equation~(6.3)]{Sim18}).
Let $\chi_{j,r}$ be the capacity cutoffs for $Z_j$ and set $Y_{j,r}:=(1-\chi_{j,r}\circ\iota_j)Y_j$.
The field $Y_{j,r}$ has compact support in $M_j$, so stability gives a nonnegative second variation for $\iota_j+tY_{j,r}$.
Since $Y_j$ and $dY_j$ are bounded, the capacity condition gives $dY_{j,r}\to dY_j$ strongly in $L^2(M_j)$ as $r\to\infty$ for each fixed $j$.
We may therefore pass to the limit and obtain the same inequality for $Y_j$.
The identity for $dY_j$, together with the oddness of $Y$ and the quadratic dependence on the differential, expresses the integrand as a fixed continuous, compactly supported function of the position and the unoriented tangent plane.
Varifold convergence allows passage to the limit in the stability inequality.
On each limiting branch, the field is $\varphi\nu$, so the normal second variation formula gives
\[
\int_{\widehat M}q|A|^2\varphi^2\,d\mathcal H_{\widehat M}^n\leq\int_{\widehat M}q|\nabla\varphi|^2\,d\mathcal H_{\widehat M}^n;
\]
see \cite[Chapter~2, \S~6]{Sim18}.
We stress that this computation uses only local unit normals and does not require $\widehat M$ to be two-sided.
Since $q$ and $|A|$ are bounded on compact subsets of $\widehat M$, intrinsic capacity cutoffs extend this inequality to every $\varphi\in C_c^\infty(\widehat M)$.

To prove stability of $M$, for $\psi\in C_c^\infty(M)$ set
\[
 h(p):=
 \left(
 \frac{1}{q(p)}
 \sum_{z\in\pi^{-1}(p)}\psi(z)^2
 \right)^{1/2}.
\]
Since $\pi:M\to\widehat M$ is a local isometry and $q$ is locally constant, $h$ is locally Lipschitz with compact support, and hence belongs to $W_c^{1,2}(\widehat M)$.
Applying the scalar inequality to $h$ by Sobolev
approximation, and using Cauchy--Schwarz for its
gradient, gives
\[
 \int_M |A_\iota|^2\psi^2\,d\mathcal H_M^n
 =
 \int_{\widehat M}q|A_{\widehat\iota}|^2h^2
 \,d\mathcal H_{\widehat M}^n
 \le
 \int_{\widehat M}q|\nabla h|^2
 \,d\mathcal H_{\widehat M}^n
 \le
 \int_M|\nabla\psi|^2\,d\mathcal H_M^n.
\]
Since $\iota$ is two-sided, this proves stability.

\medskip\noindent
\emph{Step 5: dimension bound for the singular set.}
Work in one of the graphical neighbourhoods obtained in Step~2, with limiting $q_a$-valued maps $u_a$ over the balls $B_a^n\subset P_a$.
The uniform estimate \eqref{eq:area-capacity-caccioppoli} and weak compactness give $u_a\in W^{2,2}_{\mathrm{loc}}$ for every $a$.

Fix $a$ and let $u_{a,\ell}$ be a local $C^1$ selection of $u_a$ away from $\mathcal B_{u_a}$.
We want to show that $u_{a,\ell}$ is a smooth solution of the minimal surface system.
Then $j_1u_{a,\ell}=(u_{a,\ell},Du_{a,\ell})$ is a continuous selection of $j_1u_a\in W^{1,2}_{\mathrm{loc}}$.
On almost every coordinate line, $j_1u_a$ is absolutely continuous with square-integrable derivative.
Fix such a line, parametrized by $t$, and use the same notation for the restricted maps.
For $s<t$, continuity confines $j_1u_{a,\ell}([s,t])$ to one connected component of the union of the closed balls of radius $r:=\int_s^t |(j_1u_a)'(\tau)|\,d\tau$ centred at the values of $j_1u_a(s)$.
A finite chain of these balls gives $|j_1u_{a,\ell}(t)-j_1u_{a,\ell}(s)|\le C\int_s^t |(j_1u_a)'(\tau)|\,d\tau$, where $C$ depends only on the multiplicity.
Thus $j_1u_{a,\ell}$ is absolutely continuous on almost every coordinate line, with $|\partial_i(j_1u_{a,\ell})|\le C|\partial_i(j_1u_a)|$ almost everywhere.
Consequently, $j_1u_{a,\ell}\in W^{1,2}_{\mathrm{loc}}$ and $u_{a,\ell}\in W^{2,2}_{\mathrm{loc}}$.
Lemma~\ref{lem:area-strong-equation}, applied to $u_a$, and the Sobolev chain rule then give the minimal surface equation in weak form for each selection, so elliptic regularity makes these selections smooth.
Hence, away from the lifted branch sets of the maps $u_a$, the varifold $V$ is locally a finite integer sum of smooth embedded minimal hypersurfaces.
In particular, every point of $\operatorname{sing}V$ in this neighbourhood lies over the branch set of at least one $u_a$.

By \cite[Theorem~5.7]{Luc26}, $\dim_{\mathcal H}\mathcal B_{u_a}\le n-2$ for every $a$.
Writing $u_a$ in terms of its ordered Lipschitz scalar selections shows that its graph over $\mathcal B_{u_a}$ is a finite union of Lipschitz images of $\mathcal B_{u_a}$, and therefore also has Hausdorff dimension at most $n-2$.
A countable graphical cover of $\operatorname{spt}\|V\|\setminus T$ gives $\dim_{\mathcal H}\bigl(\operatorname{sing}V\setminus T\bigr)\le n-2$.
Since $\dim_{\mathcal H}T\le n-3$ when $n\ge3$, and $T=\varnothing$ when $n=2$, we conclude that
\[
 \dim_{\mathcal H}\operatorname{sing}V\le n-2.
\]

\medskip\noindent
\emph{Step 6: conditional compactness.}
Assume Conjecture~$(\star)$ and set $Z:=\operatorname{sing}V$.
It suffices to prove $\operatorname{cap}_{2,V}(Z)=0$.

Work in one of the graphical neighbourhoods obtained in Step~2, with limiting maps $u_a$ over $B_a^n\subset P_a$.
By Step~5, every point of $Z$ in this neighbourhood lies over the branch set of at least one $u_a$.

Write each $u_a$ in terms of its ordered Lipschitz scalar selections, and use the common graphical projection constructed in Step~2.
Denote its base ball by $B^n$ and the projection by $\pi_{B^n}$.
Parametrizing each scalar graph over $B_a^n$ and applying $\pi_{B^n}$ gives a bi-Lipschitz map from the corresponding portion of $B_a^n$ onto $B^n$.
For any compact set $K$ in a smaller cylinder, Conjecture~$(\star)$, bi-Lipschitz invariance and finite subadditivity therefore give zero $2$-capacity in $B^n$ for $S:=\pi_{B^n}(Z\cap K)$.

Choose $\psi_r\in C_c^1(B^n)$ with $0\le\psi_r\le1$, $\psi_r=1$ near $S$, and $\|\psi_r\|_{W^{1,2}(B^n)}\to0$. Let $\zeta\in C_c^1(U)$ satisfy $0\le\zeta\le1$, be supported in the graphical neighbourhood, and equal one near $K$. Then $\chi_r:=\zeta\,(\psi_r\circ\pi_{B^n})$ equals one near $Z\cap K$. Since all the graphical pieces are Lipschitz over $B^n$, the area formula gives
\[
\int_U\left(\chi_r^2+|\nabla^V\chi_r|^2\right)\,d\|V\|\le C\int_{B^n}\left(\psi_r^2+|D\psi_r|^2\right)\,d\mathcal H^n\longrightarrow0.
\]
Hence $\operatorname{cap}_{2,V}(Z\cap K)=0$ for every compact $K\Subset U\setminus T$, by a finite graphical cover. By countable subadditivity, $\operatorname{cap}_{2,V}(Z\setminus T)=0$.

For $n\ge3$, monotonicity gives $\|V\|(B_r^{n+1}(x))\le C_Kr^n$ locally uniformly, while $\dim_{\mathcal H}T\le n-3$ implies $\mathcal H^{n-2}(T)=0$. Ball-covering cutoffs, with relative energy bounded by $C_K(r^{n-2}+r^n)$ for each ball, therefore give
$$
\operatorname{cap}_{2,V}(T)=0
$$
as well. For $n=2$, we have $T=\varnothing$.
Subadditivity now yields $\operatorname{cap}_{2,V}(Z)=0$, hence \eqref{eq:vanishing-capacity} and $\|V\|(Z)=0$. Thus $\iota$ represents $V$ after zero extension across $Z$, and $V\in\mathcal V(U)$, completing the proof.
\end{proof}

\begin{remark}
Let $\mathcal V_{\mathcal H^{n-2}}(U)$ denote the class
obtained by replacing condition~4 in the definition
of $\mathcal V(U)$ with
\[
\mathcal H^{n-2}(Z\cap K)<\infty
\qquad\text{for every compact }K\Subset U.
\]
This is the local formulation of the immersion class
considered in \cite{Bel25,WZ26}.
Stationarity extends across $Z$ by the removability
theorem \cite[Theorem~7.2]{Ber19}; see also
\cite[Proposition~B.5]{BLPS26}.
The capacitary cutoffs provided by
\cite[Corollary~B.6]{BLPS26} then give
\[
\mathcal V_{\mathcal H^{n-2}}(U)\subset\mathcal V(U).
\]

Compactness of this smaller class requires additional
control.
Even a theorem giving countable $(n-2)$-rectifiability
and locally finite $\mathcal H^{n-2}$-measure for the
branch set of every stationary $C^{1,\alpha}$
multiple-valued graph would provide this finiteness
only within the graphical neighbourhoods outside $T$.
For a limit $V$ of a sequence
$V_j\in\mathcal V_{\mathcal H^{n-2}}(U)$, these local
conclusions alone would not rule out accumulation
of branch-set measure near $T$, and hence would not
by themselves imply
\[
\mathcal H^{n-2}(\operatorname{sing}V\cap K)<\infty
\qquad\text{for every compact }K\Subset U.
\]
Controlling this accumulation is a further issue
of summability, which is absent from the capacitary
formulation because a countable union of sets of zero
capacity still has zero capacity.
In dimension $2$ this particular obstruction does not
arise, since $T=\varnothing$.
\end{remark}

\bigskip

\noindent
Department of Pure Mathematics and Mathematical Statistics, University of Cambridge\\
\textit{Email}: \texttt{ml2003@cam.ac.uk}

\end{document}